\documentclass[10pt]{amsart}

\usepackage{amssymb}
\usepackage[all]{xy}
\usepackage{mathrsfs}
\usepackage{enumerate}

\newcommand{\rar}[1]{\stackrel{#1}{\longrightarrow}}

\newcommand{\oX}{{\overline X}}
\newcommand{\oZ}{{\overline Z}}
\newcommand{\oK}{{\overline K}}
\newcommand{\oj}{{\overline j}}
\newcommand{\of}{{\overline f}}
\newcommand{\oC}{{\overline C}}

\newcommand{\ok}{{\overline k}}

\newcommand{\bA}{{\mathbb A}}

\newcommand{\bC}{{\mathbb C}}

\newcommand{\bF}{{\mathbb F}}
\newcommand{\bG}{{\mathbb G}}

\newcommand{\bP}{{\mathbb P}}
\newcommand{\bQ}{{\mathbb Q}}
\newcommand{\bR}{{\mathbb R}}

\newcommand{\bZ}{{\mathbb Z}}

\newcommand{\cD}{{\mathcal D}}

\newcommand{\cF}{{\mathcal F}}

\newcommand{\cH}{{\mathcal H}}

\newcommand{\cO}{{\mathcal O}}

\newcommand{\cS}{{\mathcal S}}

\newcommand{\cV}{{\mathcal V}}

\newcommand{\cX}{{\mathcal X}}

\newcommand{\sK}{{\widehat {K}}}

\newcommand{\nc}{\newcommand}

\nc\wh{\widehat}

\nc\on{\operatorname}

\nc\Gr{\on{Gr}}

\nc\Fl{\on{Fl}}

\newtheorem{conj}{Conjecture}

\newcommand{\limto}{{\displaystyle\lim_{\longrightarrow}}}
\newcommand{\rightlim}{\mathop{\limto}}

\newcommand{\leftlim}{\mathop{\displaystyle\lim_{\longleftarrow}}}
\newcommand{\limfromn}{\leftlim\limits_{\raise3pt\hbox{$n$}}}
\newcommand{\limton}{\rightlim\limits_{\raise3pt\hbox{$n$}}}

\newcommand{\rightlimit}[1]{\mathop{\lim\limits_{\longrightarrow}}\limits%
                    _{\raise3pt\hbox{$\scriptstyle #1$}}}

\newcommand{\leftlimit}[1]{\mathop{\lim\limits_{\longleftarrow}}\limits%
                    _{\raise3pt\hbox{$\scriptstyle #1$}}}

\newcommand{\epi}{\twoheadrightarrow}
\newcommand{\iso}{\buildrel{\sim}\over{\longrightarrow}}
\newcommand{\mono}{\hookrightarrow}

\DeclareMathOperator{\Coker}{{Coker}}

\DeclareMathOperator{\Ker}{{Ker}}

\makeatletter

\newcommand{\Rmnum}[1]{\expandafter\@slowromancap\romannumeral #1@}
\makeatother

\newtheorem*{theo}{Theorem}

\newtheorem{Th}{Theorem}
\newtheorem{pr}{Proposition}[section]
\newtheorem{lm}[pr]{Lemma}

\newtheorem{cor}[pr]{Corollary}

\newtheorem{example}[pr]{Example}

\theoremstyle{definition}

\newtheorem{rem}[pr]{Remark}
\newtheorem{assertion}[pr]{Assertion}

\numberwithin{equation}{section}

\newcommand{\Log}{\operatorname{Log}}

\begin{document}

\title{Motivic integral of K3 surfaces over a non-archimedean field}

\author{Allen J. Stewart and Vadim Vologodsky}



%\date{}

\begin{abstract}
We prove a formula expressing the motivic integral (\cite{ls}) of a K3 surface over $\bC((t))$ with semi-stable reduction in terms of the associated limit mixed Hodge structure. Secondly, for every smooth variety over a 
complete discrete valuation field we define an analogue of the monodromy pairing, constructed by Grothendieck  in the case of Abelian varieties, and prove that our monodromy pairing is a birational invariant of the variety. Finally, we propose a conjectural formula for the motivic integral of maximally degenerate K3 surfaces over an arbitrary complete discrete valuation field  and prove this conjecture for Kummer K3 surfaces. 
\end{abstract}

\maketitle

\section{Introduction}
%Exact cohomological  formula  for $\int_{\cX(K)}$  is not known even $X$ is an abelian variety. (See, however, \cite{} for the semi-stable case). 
   \subsection{Motivic Integral of a Calabi-Yau variety.}  Let $R$ be a  complete discrete valuation ring with fraction field $K$ and perfect residue field $k$.
  By a Calabi-Yau variety $X$ over $K$ we mean a smooth projective scheme $X$ over $K$, of pure dimension $d$, with trivial canonical bundle $\omega_X:=\Omega_{X/K}^d$. 
  In (\cite{ls}),  Loeser and Sebag  associated with any Calabi-Yau variety $X$ over $K$ a canonical element
   $$\int_{X}\in K_0(Var_k)_{loc}$$
   of the ring  $K_0(Var_k)_{loc}$, where $K_0(Var_k)_{loc}$ is obtained from the Grothendieck ring $K_0(Var_k)$ of algebraic varieties over $k$ by inverting the class $[\bA ^1_k]$ of the affine line.
   
   The motivic integral $\int_{X}$ can be computed from a weak N\'eron model of $X$.  Recall, that a weak N\'eron model of a smooth proper scheme $X$ over $K$ is a smooth scheme  $\cV$ of finite type over $R$ together with an isomorphism $\cV\otimes_R K \simeq X$ satisfying the following property: for every finite unramified extension $R'\supset R$ with fraction field $K'$, the canonical map $\cV(R') \to X(K')$ is bijective. 
   According to (\cite{blr}, \S 3.5, Theorem 3), every smooth proper $K$-scheme $X$  admits a weak N\'eron model. We note that a weak N\'eron model is almost never unique:
  for example, if $\oX$ is a proper regular model of $X$ over $R$, then the smooth locus
  $\oX_{sm}$ of $\oX$ is a weak N\'eron model of $X$ (see Lemma \ref{exampneron}). 
   
 Given a Calabi-Yau variety $X$ over $K$, a weak N\'eron model $\cV$ of $X$,  and a nonzero top degree differential form $\omega\in \Gamma(X, \omega_X)$, we can view $\omega$ as a rational section of the canonical bundle $\omega_{\cV/R}$ on $\cV$. The divisor of $\omega$ is supported on the special fiber $V^\circ$ of $\cV$. Thus, we can write
 \begin{equation}\label{def.mult}
 \mathrm{div} \, \omega = \sum _i m_i V_i^\circ,
 \end{equation}
where $V_1^\circ, \cdots, V_s^\circ$ are the irreducible components of the special fiber $V^\circ$.
  The motivic integral of  $X$ is defined by the formula\footnote{We note that our terminology and notation are different from those used by Loeser and Sebag. Notation for  $\int_{X}$
  in (\cite{ls}) is $[\oX]$. The name ``motivic integral''  is reserved in {\it loc. cit.} for a more general construction that associates with any smooth proper $K$-scheme $X$ and a top degree differential form 
  $\omega\in \Gamma(X, \omega_X)$     an element $\int_{X}\omega $ of a certain completion of the motivic ring   $K_0(Var_k)_{loc}$.}  
     \begin{equation}\label{n.m.i}
\int_{X}: =  \sum _i   [ V_i^\circ  ] ( m_i  - \min_i  m_i ). 
\end{equation}
Here, given an element $[Z]\in K_0(Var_k)_{loc}$ and an integer $n$, we write $[Z](n)$ for its Tate twist:
$$[Z](n): = [Z]\cdot[\bA^1]^{-n}.$$
A key result
  %discovered by Kontsevich and
   proven by Loeser and Sebag (\cite{ls}, Theorem 4.4.1) is that the right-hand side of  equation (\ref{n.m.i}) is independent of the choice of $\cV$ and $\omega$.
%Moreover, the motivic integral is a birational invariant of $X$. 

If $k=\bF_q$,  the image of the motivic integral under the homomorphism 
 \begin{equation}\label{p.real}
K_0(Var_{\bF_q})_{loc}\to \bZ_{(q)} \quad [Z] \rightsquigarrow |Z(\bF_q)| 
\end{equation}
is equal to the volume $\int_{X(K)}|\omega | $, for an appropriately normalized  $\omega\in \Gamma(X, \omega_X)$ (\cite{ls}, \S 4.6).

In this paper we express the motivic integral of K3 surfaces over $\bC((t))$ with strictly semi-stable reduction in terms of the associated limit mixed Hodge structures. We also compute the motivic integral
of some K3 surfaces over an arbitrary complete discrete valuation field. To our knowledge the only class of varieties, for which similar formulas were previously known, is the class of abelian varieties (see, {\it e.g.} \cite{sga7}, Expos\'e {\rm IX}, \cite{v}, \cite{hn1}, \cite{hn2}), where the computation is based on the theory of N\'eron models, and, in particular, for $K=\bC((t))$, on the Hodge theoretic description of the special fiber of the N\'eron model. 
Unfortunately,  K3 surfaces do not have a N\'eron model, in general, which makes our problem substantially more difficult.

Let us describe the organization of the paper in more detail. 
\subsection{Limit mixed Hodge structure.} In \S \ref{prel} we explain some preliminary material, the most important of which is the notion of limit mixed Hodge structure associated with a variety over the field of formal Laurent series $ \bC((t))$.
  Schmid and Steenbrink  associated with every smooth projective variety over the field $K_{mer}$ of meromorphic functions  on an open neighborhood of zero in the complex plane a mixed Hodge structure, called the limit mixed Hodge structure. In \S \ref{l.h.s},  using Log Geometry,  we extend  the Steenbrink-Schmid construction to smooth projective varieties over  $\bC((t))$.  
 
   \subsection{Motivic integral of K3  surfaces over $\bC((t))$.}\label{MR}
   In order to state our first main result we need to introduce a bit of notation. Let $X$ be a smooth projective K3 surface over $K=\bC((t))$ and let 
   %some arithmetic invariants of the monodromy operator.   
%Let $\overline C$ be a smooth curve over $\bC$, $a\in \overline C$ a point, $C=\overline C -a$,  $t$ a local coordinate at $a$, and let  $X\rar{\pi} C$ be smooth projective scheme over $C$ of relative dimension $n$.  
 $H^2(\lim X)= (H^2(\lim X,\bZ), W_i^{\bQ}, F^i)$
 be the corresponding limit mixed Hodge structure (see  \S \ref{l.h.s}).  Assume that the monodromy acts on $H^2(\lim X,\bZ)$ by a unipotent operator.  Then its logarithm $N$ is known to be integral (\cite{fs}, Prop. 1.2):  
  \begin{equation}\label{logmonod}
 N:  H^2(\lim X, \bZ)  \to H^2(\lim X, \bZ).
\end{equation}
 Set $ W_i^{\bZ}= W_i^{\bQ}\cap H^2(\lim X,\bZ)$. The morphisms 
  \begin{equation}\label{m.p}
Gr \, N^i:   W_{i+2}^{\bZ}/W_{i+1}^{\bZ}  \to  W_{2-i}^{\bZ}/W_{1-i}^{\bZ}, \quad i=1, 2 
\end{equation}
are injective and have finite cokernels.  Let $r_{i}(X, K)$ be their orders. In \S \ref{p.m.t} we prove the following result.
\begin{Th}\label{main}
Let $X$ be a smooth projective K3 surface over $K=\bC((t))$.  Assume that $X$  has a strictly semi-stable model over $R=\bC[[t]]$ and that the operator $N$ is not equal  to $0$.
 Let $s$ be the smallest integer such that $N^s= 0$. Then $s$ is either $2$ or $3$ and for every finite extension $K_e\supset K$ of degree $e$ the motivic integral of the K3 surface  $X_e= X \otimes_K K_e$ over $K_e$ is given by the following formulas.
\begin{enumerate}[(a)]
\item If $s=2$ then
  \begin{eqnarray}\label{eq1m}
\int_{X_e}= 2\bZ(0) -(e\sqrt{r_{1}(X, K)}+1)[E(X)]+20\bZ(-1)\qquad \qquad &  \\
\nonumber  \qquad\qquad\qquad +(e\sqrt{r_{1}(X, K)}-1)[E(X)](-1)+2\bZ(-2),&
  \end{eqnarray}
where $E(X)$ is the elliptic curve defined by the weight  $1$ Hodge structure on $W_1^{\bZ}=  W_1^{\bQ}\cap H^2(\lim X,\bZ)$ and $\bZ(n): = [\bA^1]^{-n}$, $n\in \bZ$.
\item If $s=3$ then 
  \begin{eqnarray}\label{eq2m}
\int_{X_e}= \left(\frac{e^2 r_{2}(X, K)}{2}+2\right)\bZ(0)+(20-e^2 r_{2}(X, K))\bZ(-1)\qquad \quad &\\
\nonumber +\left(\frac{e^2 r_{2}(X, K)}{2}+2\right)\bZ(-2).&
  \end{eqnarray}
 \end{enumerate}
  \end{Th} 
Note, that if $N=0$ the K3 surface $X$  has a smooth proper model over $R$    
 whose special fiber $Y$ (and thus the motivic integral) is determined by the polarized pure Hodge structure $H^2(\lim X,\bZ)$.
%\end{rem}
%\begin{rem} In \S \ref{a.v.k.s} we prove, that for every family of Kummer surfaces $W_{0}^{\bZ}=  W_0^{\bQ}\cap H^2(\lim X,\bZ)$.  We do not know if this is true for an arbitrary family of K3 surfaces.
%\end{rem} 

Let us explain the idea of our proof assuming that $e=1$.   
%Let $\overline C$ be a smooth curve over $\bC$, $a\in \overline C$ a point, $C=\overline C -a$, and let $t$ be a local coordinate at $a$. 
 First, using the theory of Hilbert schemes and Artin's approximation theorem, we reduce the proof to the case when $X$ is obtained by the restriction of a smooth family $\cX$ of K3 surfaces over a smooth punctured complex curve $C=\overline C -a $ to the formal punctured neighborhood of
 the point $a\in \overline C$.
 The rest of the proof is based on a result of Kulikov (\cite{ku}) asserting the existence of a (non-unique) strictly semi-stable model  $\overline{\cX}\rar{\overline \pi} \overline{C}$ such that the log canonical bundle $\omega _{\overline{\cX}/\overline C}(log)$ is trivial over an open neighborhood of 
the special fiber $Y$. For any such model, we have
$$\int_{X} =   [ Y_{sm}  ], $$
where $Y_{sm} \subset Y$ is the smooth locus of $Y$.
 It is shown in (\cite{ku}) that the special fiber $Y$ of a Kulikov model has a very special form. 
 If $s=2$ the Clemens polytope $Cl(Y)$ of $Y$ (see \S \ref{c.p.}) is  a partition of an interval and all but two irreducible components of $Y$ are ruled surfaces fibered over elliptic curves, all of which are isomorphic to a single elliptic curve $E$.   
 The two components corresponding to the boundary points of $Cl(Y)$
 are rational surfaces. If $s=3$ then all the irreducible components of $Y$ are rational surfaces and the Clemens polytope $Cl(Y)$ is a triangulation of a sphere. Next, using results of  Friedman and Scattone (\cite{fs}, \cite{fr}) we prove that
  the Steenbrink weight spectral sequence for Kulikov's  model $\overline{\cX}\rar{\overline \pi} \overline{C}$ (and therefore by the Weak Factorization Theorem (\cite{konso}, Theorem 9), for every strictly semi-stable model of $X$) degenerates integrally at the second term. Of course, the degeneration of the weight spectral sequence with rational coefficients is a corollary of Hodge Theory and holds in general, but the degeneration over $\bZ$ is a special non-trivial property of  K3 surfaces. 
  This, combined with the generalized Picard-Lefschetz formula, implies that, for $s=2$,  the Hodge structure on $H^1(E)$ is isomorphic to that on $W_1^{\bZ}$ and that the number of irreducible components of $Y$ equals $\sqrt{r_1(X, K)} +1$.
  Similarly, as proven in  (\cite{fs}), for $s=3$, the combinatorics of $Y$ ({\it i.e.}, the number of irreducible, components, double curves and triple points) is completely determined by the monodromy action on the integral lattice $H^2(\lim X,\bZ)$.
  This, together with a variant of  A'Campo's formula (Proposition \ref{modQ(1)-Q}), completes the proof. 
  
     %\footnote{Thanks to the Torelli type results for this class of varieties, we can compute even the motivic integral itself. } We also propose a conjectural formula for the motivic integral of semi-stable K3 surfaces over an arbitrary complete discrete valuation domain and prove this formula for %Kummer K3 surfaces. 
%We also propose a conjectural formula for the volume of semi-stable K3 surfaces over local fields. We begin with the story over $\bC$.  

\subsection{Monodromy pairing.}\label{mpi}  
In \S \ref{s.m.p.} we introduce a generalization of the invariant $r_2(X,K)$, that we defined in \S \ref{MR} for K3 surfaces over $\bC((t))$, to the case of an arbitrary smooth variety over a complete discrete valuation field.  Our construction is based on the theory of analytic spaces  over  non-archimedean fields developed
by Berkovich (\cite{ber1}).  For a complete discrete valuation field $K$ we denote by $\sK$ the completion of an algebraic closure of $K$.
 One of the key features of Berkovich's theory is that the underlying topological space $|X^{an}_\sK |$ of the analytification of a scheme $X$ over  $K$ has interesting topological invariants (in contrast with the space $X(\sK)$ equipped with the usual 
 topology, which is totally disconnected).  In particular, if $X$ is the generic fiber of a proper strictly semi-stable scheme $\oX$ over $R$ the space $|X^{an}_\sK |$ is homotopy equivalent to the Clemens polytope of the special fiber $Y$. 
 We denote by $\Gamma^m_C (X)$ the singular cohomology of the space $|X^{an}_\sK |$ with coefficients in a ring $C$.  In Theorem \ref{main2}, we prove that, for every prime $\ell$ different from the characteristic of the residue field of $K$, and for every smooth scheme $X$,
  there is a canonical isomorphism of $\mathrm{Gal}(\oK/K)$-modules
 \begin{equation}\label{mibi}
\gamma: \Gamma^m_{\bQ_{\ell}}(X) \iso  Im\left( H^m( X_{\overline K}, \bQ_{\ell})(m)\rar{N^m} H^m( X_{\overline K}, \bQ_{\ell})\right),
\end{equation}
where $N$ is the logarithm of the monodromy operator.  In particular, the dimension of the vector space on the right-hand side of  (\ref{mibi}) is independent of $\ell$. Let us note, that a different description of the space $\Gamma^m_{\bQ_{\ell}}(X) $ in the case of finite residue field
was obtained earlier by Berkovich (\cite{ber4}). 
%The proof of (\ref{mibi}) is long and difficult.
%We first prove the theorem assuming that $X$ is projective and has a strictly semi-stable model  $\oX$ over $R$. In this case, (\ref{mibi}) is equivalent to a special case of Deligne's monodromy conjecture that we prove using the method of Steenbrink 
  %(who proved it for all weights and $k=\bC$). To prove the theorem for arbitrary smooth $X$ we show that the functors $\Gamma^m_{\bQ_l}$,  $N^m H^m$, firstly, admit transfers for finite morphisms
%and, secondly, take every dominant open embedding $U\mono X$ to an isomorphism.  Finally, we use de Jong's alteration result to complete the proof.

 If $d$ is the dimension of $X$, we use (\ref{mibi}) to define a non-degenerate  pairing
 \begin{equation}\label{monpairi}
\Gamma^d_\bQ (X) \otimes \Gamma^d_\bQ(X)\to \bQ. 
\end{equation}
In the special case when $X$ is proper, the pairing (\ref{mibi})
is given by the formula
\begin{equation}\label{monparintro}
 (x, y ) = (-1)^{\frac{d(d-1)}{2}}  <\gamma(x), y'>,
 \end{equation}
 where $y' \in  H^d( X_{\overline K}, \bQ_{\ell}) $ is an element such that $N^dy'=\gamma(y)$ and $<,>$
 is the Poincar\'e pairing on $H^d( X_{\overline K}, \bQ_{\ell}) $. We prove in Theorem \ref{main3} that (\ref{monparintro}) is independent of $\ell$ and positive.
Moreover, the groups $\Gamma^m_C (X)$ and the monodromy pairing (\ref{monpairi}) are birational invariants of $X$.

We define a numeric (birational) invariant $r_d(X,K)$ of $X$ to be the discriminant of  the dual pairing 
\begin{equation}\label{mibi2}
\Gamma_d(X) \otimes \Gamma_d(X) \to \bQ,
\end{equation}
where $\Gamma_d(X)$ is $\mathrm{Hom}(\Gamma^d_\bZ(X), \bZ)$.

 In remark \ref{m.p.p.v.},  we define for a polarized projective variety $X$ and any integer $m$ a more general positive pairing  $\Gamma^m_\bQ(X) \otimes \Gamma^m_\bQ(X) \to \bQ$ which in the case of semi-stable abelian variety $A$ and its dual $A'$ boils down, after some identifications,
 to the monodromy pairing $\Gamma_1(A) \otimes \Gamma_1(A') \to \bZ $ defined by Grothendieck (\cite{sga7}, Exp. {\rm IX}). In particular, the number $r_d(A,K)$ is non-zero if and only if $A$ is completely degenerate in which case  $r_d(A,K)$ is equal to $d!|\pi_0(\cV(A)\otimes \ok)|$, where
 $\cV(A)$ is the N\'eron model of $A$.

\subsection{Motivic integral of maximally degenerate K3 surfaces.}  We say that a $d$-dimensional Calabi-Yau variety over a complete discrete valuation field $K$ is maximally degenerate if $\Gamma^d_\bQ(X) \ne 0$. 
According to (\ref{mibi}), $X$ is maximally degenerate if and only if for some (and, hence, for any) prime $\ell\ne \mathrm{char} \, k$ the map
  $$H^d( X_{\overline K}, \bQ_{\ell})(m)\rar{N^d} H^d( X_{\overline K}, \bQ_{\ell})$$
  is not zero\footnote{There is an extensive literature on maximally degenerate Calabi-Yau varieties over $\bC((t))$. See  {\it e.g.}  \cite{mo1}, \cite{lty}.}.
 We conjecture that for every maximally degenerate K3 surface over $K$ there exists a finite extension $K'\supset K$ such that, for every finite extension $L\supset K$ of ramification index $e$ containing $K'$, we have
$$\int_{X_L}= \left(\frac{e^2 r_{2}(X,K)}{2}+2\right)\bQ(0)+(20-e^2 r_{2}(X,K))\bQ(-1)+\left(\frac{e^2 r_{2}(X,K)}{2}+2\right)\bQ(-2).$$
 If  $\mathrm{char}\, k =0$ our conjecture follows from part (b) of Theorem \ref{main}. 
 In \S  \ref{m.i.m.g.k3}     we prove this conjecture in the case of Kummer K3 surfaces over an arbitrary complete discrete valuation field $K$ with $\mathrm{char}\,  k \ne 2$ by constructing explicitly a poly-stable formal model of the analytic space
$X^{an}$. 

The groups $\Gamma^d_\bZ(X)$ that we used to define the invariant $r_d(X,K)$ can be interpreted as the weight $0$ part of the limit motive of $X$ (Remark \ref{w.0.a}). It would be interesting to define geometrically the limit 1-motive attached to $X$ and use it to compute the motivic integral for K3 surfaces which are not maximally degenerate. 
    
%    \subsection{Organization of the paper.}  In \S \ref{m.i.c.y.v} we review the construction of motivic integrals and prove two new results. First, we show that for every  smooth projective family of Calabi-Yau varieties $X\rar{\pi} C$ admitting a semi-stable model over $ \overline C$ the image of %$R^{Hodge}(\int_{X})$ in 
   % the quotient ring $ K_0(MHS)/(\bQ(1)-\bQ)$ equals to $\sum (-1)^i [H^i(\lim X)]$. Secondly, we prove that for every Calabi-Yau variety  over a local field $K$ admitting a semi-stable model over $R$ the formula (\ref{g.l.f}) holds modulo $q-1$.
%  In \S \ref{p.m.t}  we prove Theorem \ref{main}. In \S \ref{a.v.k.s} we compute the motivic integral for arbitrary abelian varieties and K3 surfaces over $\bC((t))$. Finally, in the last paragraph we propose a conjectural formula for the volume of semi-stable K3 surfaces over  local fields. 

 \section*{Acknowledgements}
 We are grateful to David Kazhdan, who asked the second author to write a cohomological formula for the p-adic measure of a Calabi-Yau variety over $\bQ_p$ and suggested to work out the case of K3 surfaces, and to Vladimir Berkovich
for answering our numerous questions on Non-archimedean Analytic Geometry and his help with \S \ref{s.m.p.}.  Special thanks go to the referee for careful reading the first draft of the paper
and for his (or her)  numerous remarks and suggestions.

 \section{Preliminaries.}\label{prel}
  \subsection{Clemens Polytope and nerve of a strictly semi-stable scheme}\label{c.p.} 
  %{\bf Add: poly-stable schemes,  Berkovich's result, Saito's Lemma 1.11}
   Let $R$ be a complete discrete valuation ring with residue field $k$ and fraction field $K$. 
   Recall that a scheme $\oX$ of finite type over $\mathrm{spec}\ R$ is strictly semi-stable if every point $x\in \oX$ has a Zariski neighborhood $x\in U \subset \oX$ such that the morphism
   $U \to \mathrm{spec}\, R$ factors through an \'etale morphism 
   $$U\to \mathrm{spec} \, R[T_0,\ldots, T_{d}]/(T_0\cdots T_{r}-t), \quad 0\leq r \leq d, $$
       for a uniformizer  $t$ of $K$.  
       %If, for every point $x$, the number $s$ in the above definition is $1$ the scheme $\oX$ is called strictly semi-stable over $R$.
     %Let $R$ be a complete discrete valuation ring with residue field $k$ and fraction field $K$. 
  % We say that a scheme $\oX$ of finite type over $\mathrm{spec}\ R$ is strictly poly-stable if every point $x\in \oX$ has a Zariksi neighborhood $x\in U \subset \oX$ such that the morphism
   %$U \to \mathrm{spec}\, R$ factors through an \'etale morphism $U\to \mathrm{spec} \, A_1 \times_R \ldots \times_R \mathrm{spec} \, A_s$, where each  $A_i$ is of the form
   %$$R[T_0,\ldots, T_{n_i}]/(T_0\cdots T_{r_i}-t), \quad 0\leq r_i\leq n_i$$
    %for a uniformizer  $t$ of $K$.  If, for every point $x$, the number $s$ in the above definition is $1$ the scheme $\oX$ is called strictly semi-stable over $R$.
     If $k$ is perfect, $\oX$ is a strictly semi-stable scheme if and only if
it is regular and flat over $R$, the generic fiber $X=\oX\times _R K$ is smooth over $K$ and the special fiber $Y=\oX \times _R k$ is a reduced strictly normal crossing divisor on $\oX$.
%Similarly, one defines the notion of poly-stable (resp. semi-stable) formal scheme over $R$ ({\it e.g.}, \cite{ber3},  Def. 1.2).
  
  Let $\oX$ be a  strictly semi-stable scheme.  Then the irreducible components  $V_1, \ldots, V_m$ of $Y$ as well as the schemes 
    \begin{equation}\label{sstm}
Y^{(q)}=\coprod_{i_0<\cdots<i_q} V_{i_0}\cap\cdots\cap V_{i_q}
 \end{equation}
  are smooth. It is convenient to encode the combinatorial structure of $Y$ by a certain topological space.  To do this we need to introduce some terminology.
   
  By an abstract triangulated set we mean a contravariant functor $\tilde \Delta \to \text{Sets} $,  where $\tilde \Delta$ is the category whose objects are finite totally ordered 
  sets $[q]:=\{0, \cdots, q\}$, $q\in \bZ_{\ge 0}$ and whose morphisms are strictly increasing maps. Thus, giving an abstract triangulated set $S_{\bullet}$ amounts to giving
  a set $S_q$ of ``$q$-simplices'' for each $q\in \bZ_{\ge 0}$ together with ``boundary maps''  $\delta_j:   S_q \to S_{q-1}$, $j=0,\cdots, q$, subject to certain conditions\footnote{The category of abstract triangulated sets can be viewed as a full subcategory of the category of simplicial sets: if $S'_{\bullet}$ is a simplicial set such that the boundary of each nondegenerate simplex of $S'_{\bullet}$ is nondegenerate then nondegenerate simplices of $S'_{\bullet}$ together with the boundary maps form an abstract triangulated set.
  This yields an equivalence between the full subcategory of the category of simplicial sets whose objects satisfy the above property and the category of abstract triangulated sets (\cite{gm}, \S 1.6).}.  We shall write $|S_{\bullet}|$ for the realization of $S_{\bullet}$ (\cite{gm}, \S 1.1).
    
 Given a  strictly semi-stable scheme $\oX$ consider the abstract triangulated set whose $q$-dimensional simplices are indexed by the set $\pi_0\left(Y^{(q)}_\ok\right)$.  The boundary maps  $\delta_j: \pi_0\left(Y^{(q)}_\ok\right)\to \pi_0\left(Y^{(q-1)}_\ok\right)$, $j=0,\ldots, q$,
are given by the maps
$$\pi_0\left(V_{i_0, \ok}\cap\cdots\cap V_{i_q, \ok}\right)\to \pi_0\left(V_{i_0, \ok}\cap\cdots\cap V_{i_{j-1}, \ok}\cap V_{i_{j+1}, \ok}\cap \cdots \cap V_{i_q, \ok}\right)$$
induced by the injections 
$$V_{i_0}\cap\cdots\cap V_{i_q}\hookrightarrow V_{i_0}\cap\cdots\cap V_{i_{j-1}}\cap V_{i_{j+1}}\cap \cdots \cap V_{i_q}.$$ 
The realization of this triangulated set  is a topological space which we call (following \cite{konso}) the  Clemens polytope of $Y$ and denote  by $Cl(Y)$. Although the abstract triangulated set we constructed depends upon the choice of ordering on the set of  
irreducible components $V_i$,  the homeomorphism type of the topological space $Cl(Y)$  does not.
\begin{pr}
Let $\oX$ be a strictly semi-stable model of $X$ over $\mathrm{spec}\ R$ with special fiber $Y$ then for every abelian group $C$,
$$H^*_{sing}(Cl(Y), C)\cong H^*_{Zar}(Y_\ok,\underline{C})$$
\end{pr}
\begin{proof}  To simplify our notation we assume that $k=\ok$. Consider the complex
$$ {i_0}_*\underline{C}\rar{\partial_0} {i_1}_*\underline{C}\rar{\partial_1}\cdots.$$
where $i_q:Y^{(q)}\hookrightarrow Y$. The differentials $\partial_i$ are characterized by the property that the induced map on global sections $\Gamma({i_q}_*\underline{C})=C[\pi_0(Y^{(q)})]\to C[\pi_0(Y^{(q+1)})]= 
\Gamma({i_{q+1}}_*\underline{C})$ equals $\sum_{j}(-1)^j\delta_{j}^*$. This complex is a resolution of $\underline{C}$. Since each $Y^{(q)}$ is a
 disjoint union of smooth irreducibles and constant sheaves on irreducibles are flabby it follows that the sheaves ${i_q}_*\underline{C}$ are flabby. Thus the complex of global sections 
$$ \Gamma({i_0}_*\underline{C})\rar{\partial_0}\Gamma( {i_1}_*\underline{C})\rar{\partial_1}\cdots.$$
computes the Zariski cohomology  $H^*_{Zar}(Y,\underline{C})$.  On the other hand, this complex is the simplicial  complex of $Cl(Y).$ 

\end{proof}
   
Assume that $\oX$ is a proper semi-stable scheme over $R$. Then, by the Proper Base Change theorem ({\it e.g.}, \cite{d}), for every torsion abelian group $C$ we have  canonical morphisms of $\mathrm{Gal}(\oK/K)$-modules
 \begin{equation}\label{eq666}
  H^*_{Zar}(Y_\ok,\underline{C})\to H^*_{et}(Y_\ok,\underline{C})\cong H^*_{et}(\oX_{R^{sh}},\underline{C})\to H^*_{et}(X_\oK,\underline{C}),
  \end{equation}
  where  $R^{sh}$ denotes a strict Heselization of $R$. Applying (\ref{eq666}) to $C=\bZ/{\ell}^n\bZ$ and passing to the limit, we obtain a canonical morphism 
    \begin{equation}\label{mclet}
  H^*_{sing}(Cl(Y), \bZ_{\ell}) \to  H^*(X_\oK,\bZ_{\ell}).
 \end{equation}
 We will see in \S \ref{s.m.p.} that the groups $H^*_{sing}(Cl(Y), C)$ and the morphism (\ref{mclet}) depend only on the generic fiber $X$ and not on the choice of proper strictly semi-stable model $\oX$.
\begin{rem}\label{nerve}
Let us explain the relation of the notion of Clemens polytope to a more general notion of nerve of a scheme, introduced in (\cite{ber3}).
 For a reduced scheme $Y$ over $k$,  let $\mathrm{Nor}(Y)\subset Y $ be the normal locus of $Y$, which is an open subset of $Y$, and let $Y^{[0]}=Y$,  $Y^{[i+1]}=Y^{[i]}\backslash \mathrm{Nor}(Y^{[i]})$, $i\geq 0$.
The irreducible components of  $Y^{[i]}\backslash Y^{[i+1]}$ are called strata of $Y$.  The set, $\mathrm{Str}(Y)$, of all strata has a natural partial order:  for strata $x,y\in \mathrm{Str}(Y)$, we say that $x\leq y$ if $y$ is contained in the closure of $x$. We denote by $N(Y)$ the nerve
of the partially ordered set $\mathrm{Str}(Y)$. If $\oX$ is a strictly semi-stable scheme over $R$, the triangulated space $|N(Y\otimes \overline k)|$ is obtained from $Cl(Y)$ by subdivision. In particular, the spaces   $|N(Y\otimes \overline k)|$ are $Cl(Y)$ homeomorphic.
\end{rem}

% and that (\ref{mclet})  is injective modulo torsion
%for $l\ne \mathrm{char} \, k$.
 % In the case where $R$ is not complete we say that a scheme $\oX$ of finite type over $\mathrm{spec}\ R$ is strictly semistable if $\oX$ is strictly semistable over $\mathrm{spec}\ \hat{R}$, where $\hat{R}$ is the completion of $R$.
% If the residue field $k$ is perfect then the converse is true.
%The goal of this section is to prove the following theorem
%\begin{Th}
%Let $X$ be a strictly semi-stable Noetherian scheme over $K$ with special fiber $Y$ then for any Abelian group $A$ there exists a topological space, $Cl(Y)$, such that 
%$$H^*_{sing}(Cl(Y),A)\cong H^*_{Zar}(Y,\underline{A})$$
%\end{Th}  

\subsection{The limit mixed Hodge structure associated with a variety over $\bC((t))$.}\label{l.h.s}
In  (\cite{st1}),  Steenbrink  associated with every smooth projective variety over the field $K_{mer}$ of meromorphic functions  on an open neighborhood of zero in the complex plane a mixed Hodge structure, called the limit mixed Hodge structure. Another construction of the same
mixed Hodge structure had been given earlier by Schmid  (\cite{sch}). In this section, we explain how to extend  the Steenbrink-Schmid construction to smooth projective varieties over the field of formal Laurent series $K=\bC((t))$. 
A rough idea: generalizing a construction by Steenbrink (\cite{st3}) we attach a mixed Hodge structure to every projective normal crossing (not necessarily reduced)  log scheme over the log point.  Applying this construction to the special fiber $Y$ of a normal crossing model $\oX$ of $X$ over $R= \bC[[t]]$ we get our   $H^m(\lim X)$. We then prove independence of the choice of a model and functoriality.

We shall summarize the properties  of our construction
in the following theorem.
\begin{Th}\label{lhs}
For every non-negative integer $m$, there exists a contravariant functor
  \begin{equation}\label{fhs}
SmPr_K \to \widetilde {MHS}, 
\end{equation}
$$X \rightsquigarrow H^m(\lim X)= (H^m(\lim X,\bZ), W_i^{\bQ}, F_i, T)$$
from the category of smooth projective varieties over $K=\bC((t))$ to the category of mixed Hodge structures equipped with an endomorphism $T$ of the underlying abelian group with the following properties.
\begin{enumerate}[(a)]
\item
 If  we write $ T_{\bQ}= SU$ for the factorization  of the endomorphism $T_\bQ \in {\mathrm {End}}( H^m(\lim X,\bQ))  $ into the product of  semi-simple and unipotent  endomorphisms,  $S$ and $U$ respectively, 
 such that  $ST_\bQ = T_\bQ S$ and $UT_\bQ=T_\bQ U$,   
 then $N=log U$  is a morphism of rational mixed Hodge structures  
 $$N:  H^m(\lim X) \otimes \bQ \to H^m(\lim X) \otimes \bQ(-1)$$
 and $S$ is a finite order automorphism of   $ H^m(\lim X) \otimes \bQ $.
%  \item If, in addition,  $X$ is projective the morphism
% $$  Gr \, N^i:   W_{i+m}^{\bQ}/W_{i+m-1}^{\bQ}  \to  W_{m-i}^{\bQ}/W_{m-i-1}^{\bQ}$$
% is an isomorphism.
%\item
%The functor  (\ref{fhs}) is compatible with the base change. That is,  if $K_e= \bC((t^{\frac{1}{e}}))\supset K$ is a finite extension and  $X \rightsquigarrow  X_{K_e}$ is the base change functor, we have a functorial isomorphism
%$$(H^m(\lim X,\bZ), W_i^{\bQ}, F^i, T^e) \simeq (H^m(\lim X_{K_e}, \bZ), W_i^{\bQ}, F^i, T).$$
\item The functor  (\ref{fhs}) is compatible with base change. That is,  if $K_e= \bC((t^{\frac{1}{e}}))\supset K$ is a finite extension and  $X \rightsquigarrow  X_{K_e}$ is the base change functor, we have a functorial isomorphism
$$(H^m(\lim X,\bZ), W_i^{\bQ}, F^i, T^e) \simeq (H^m(\lim X_{K_e}, \bZ), W_i^{\bQ}, F^i, T).$$
\item If $\oX$ is a strictly semi-stable scheme over $R= \bC[[t]]$, $X$ and $Y$ are the generic and special fibers of $\oX$ respectively, and $Y^{(q)}\mono Y$ is the closed subscheme defined in (\ref{sstm}),
one has the weight spectral sequence $E_r^{pq}(\oX)$ which converges to $H^*(\lim X)$ in the category of mixed $\bZ$-Hodge structures with the first term given by the formula:
$$E_1^{pq}(\oX)= \bigoplus_{i, i-p \geq 0} H^{q+2p-2i}(Y^{(2i-p)})(p-i).$$
The sequence $E_r^{pq}(\oX) \otimes \bQ$ degenerates at $E_2$ terms.   
\item If $\cX$ is a  smooth projective variety over $K_{mer}$  the limit mixed Hodge structure   $H^m(\lim (\cX\otimes _{K_{mer}} K))$ is canonically isomorphic to the one constructed by Schmid and Steenbrink (\cite{sch}, \cite{st1}, \cite{st2}).
\end{enumerate}
\end{Th}
\begin{proof} Let    
$$(\oX, M_{\oX}) \to (\mathrm{spec}\, R, M_R=R-0)$$
be a proper smooth morphism of fine and saturated (fs for short) log schemes (\cite{il2}, \S 1). Assume that the log structure on  $(\oX, M_{\oX})$ is vertical {\it i.e.}, the induced log structure on $j: X=\oX \otimes _R K\mono \oX$ is trivial.
 A basic example of this situation is a regular proper $R$-scheme $\oX$ such that its reduced special fiber $Y_{red}$ is a normal crossing divisor on $\oX$ endowed with the log structure
  \begin{equation}\label{can.log.str}
  M_{\oX}= j_* \cO^*_X \cap \cO_{\oX}.
  \end{equation}

The special fiber $Y=\oX \otimes _{R} \bC$ with the induced log structure is a proper smooth log scheme over the log point
 $$\pi: (Y,M_Y) \to (\mathrm{spec} \, \bC)_{log}.$$  
 %let $Y_{an}=Y(\bC)$ be the analytic space associated with $Y$.
 % Consider the log structure $M_{\oX}$ on $\oX$ given by the closed subscheme $Y\mono X$ {\it i.e.},  $M_{\oX}= j_* \cO^*_X \cap \cO_{\oX}$, where
  %$j: X=\oX \otimes _R K \mono \oX$ is the open embedding. 
  Following (\cite{kn} \S 1) we consider the associated map of topological spaces
 $$\pi: Y^{log} \to (\mathrm{spec} \, \bC)^{log}=S^1,$$
 where $S^1\subset \bC$ is the unit circle. The map $\pi$ is a locally trivial fibration over $S^1$ (\cite{no}, Theorem 5.1).  Let $\exp(2\pi i \tau):  \bR^1\to  S^1$ be the universal cover, and let $\tilde Y^{log}$ be the fiber product
 $Y^{log}\times_ {S^1} \bR^1$. The topological space  $\tilde Y^{log}$ carries a canonical automorphism that takes a point  $(y, a)\in Y^{log}\times_ {S^1} \bR^1$ to $(y, a+2\pi i)$.  We will write $T_Y$ for the induced automorphism
 of the cohomology group    $H^m(\tilde Y^{log}, \bZ)$.
 %let  $\tilde p: \tilde Y^{log} =Y^{log}\times_ {S^1} \bR^1 \to Y_{an}$    be the composition of the projection to the first factor and the canonical map  $p: Y^{log}\to Y_{an}$.
 %The complex $R\tilde p_* \bZ \in D^b(Sh(Y_{an})) $ has a canonical automorphism $T_Y$ induced  by the automorphism of the space   $Y^{log}\times_ {S^1} \bR^1$, $(y, a)\rightsquigarrow (y, a+2\pi i)$. The complexes $R\tilde p_* \bZ $ and $R p_* \bZ$ are related by the exact triangle 
  %\begin{equation}\label{extrian}
%R p_* \bZ \rar{\alpha}  R\tilde p_* \bZ \rar{T_Y-1} R\tilde p_* \bZ\rar{\beta} R p_* \bZ[1].
%\end{equation}
 %Abusing notation, we will write $T_Y$ for the automorphism
 %of the group    $H^m(\tilde Y^{log}, \bZ) \simeq H^m(Y_{an}, R\tilde p_* \bZ) $ induced by the automorphism of the complex $R\tilde p_* \bZ$.
The following lemma implies that the cohomology of $\tilde Y^{log}$ depends only on the generic fiber of $\oX$.
  \begin{lm}\label{indep} Let $f: (\oX, M_{\oX})  \to (\oX', M_{\oX'})$ be a log morphism of  smooth proper  vertical fs log schemes over $(\mathrm{spec}\, R,  M_R)$. 
  Assume that the induced morphism of the generic fibers $f_K: \oX \otimes K \to \oX' \otimes K$ is an isomorphism. Then, for every non-negative integer $m$,  the morphism
   $$f^*: H^m(\tilde Y^{\prime log}, \bZ) \to  H^m(\tilde Y^{log}, \bZ)$$
   is an isomorphism. 
  \end{lm}
 \begin{proof} Let $n$ be a positive integer. The comparison theorems of Kato and Nakayama (see, e.g. \cite{il2}, Th. 5.9,  Cor. 8.4) imply the existence of the commutative diagram below.
      \begin{equation}
\def\normalbaselines{\baselineskip20pt
\lineskip3pt  \lineskiplimit3pt}
\def\mapright#1{\smash{
\mathop{\to}\limits^{#1}}}
\def\mapdown#1{\Big\downarrow\rlap
{$\vcenter{\hbox{$\scriptstyle#1$}}$}}
\begin{matrix}
 H^m(\tilde Y^{\prime log}, \bZ/n\bZ)      &      \iso{}    &   H^m_{et}(\oX' \otimes \oK, \bZ/n\bZ)  \cr
  \mapdown{f^* }& & \mapdown{f^*_{\oK}}  \cr
  H^m(\tilde Y^{log}, \bZ/n\bZ)   &      \iso{}  & H^m_{et}(\oX \otimes \oK, \bZ/n\bZ)       
\end{matrix}
 \end{equation}
 Since the groups $H^m(\tilde Y^{\prime log}, \bZ)$, $H^m(\tilde Y^{log}, \bZ)$ are finitely generated the lemma follows.
\end{proof}

 Let us explain how the formation $(T_Y, H^m(\tilde Y^{log}, \bZ)  )$ is compatible with base change.  For a positive integer $e$,  the fs log scheme   $(\mathrm{spec} \, R_e= \mathrm{spec} \, \bC[[t^{\frac{1}{e}}]],  M_{R_e}= R_e -0)$
is smooth over $(\mathrm{spec}\, R, M_R)$.  Let  $(\oX_e, M_{\oX_e})$ be the fiber product   
$$(\oX, M_\oX) \otimes _ {(\mathrm{spec}\, R, M_R)} (\mathrm{spec}\, R_e, M_{R_e})$$
in the category of fs log schemes\footnote{Warning: the functor that takes a fs log scheme to the underlying scheme does not commute with the fiber products.}.   As the functor $(Y, M_Y) \rightsquigarrow Y^{log}$ commutes with fiber products  we have 
 a  Cartesian diagram of topological spaces
       \begin{equation}\label{dcd10}
\def\normalbaselines{\baselineskip20pt
\lineskip3pt  \lineskiplimit3pt}
\def\mapright#1{\smash{
\mathop{\to}\limits^{#1}}}
\def\mapdown#1{\Big\downarrow\rlap
{$\vcenter{\hbox{$\scriptstyle#1$}}$}}
\begin{matrix}
   Y^{log}_e   &      \rar{}    &    Y^{log}  \cr
  \mapdown{ }& & \mapdown{}  \cr
  S^1  &      \rar{}  & S^1,      
\end{matrix}
 \end{equation}
where the lower horizontal map is an  $e$-fold cover.
 We get from (\ref{dcd10}) a canonical isomorphism   
 \begin{equation}\label{bc3}
H^m(\tilde Y^{log}, \bZ)\iso H^m(\tilde Y_e^{log}, \bZ)
 \end{equation}
 that carries $T_{Y_e}$ to $T_Y^e$.

 Assume, in addition, that the log scheme $\pi: (Y,M_Y) \to (\mathrm{spec} \, \bC)_{log}$  satisfies the following condition:\\ 
 
 (U):  for every closed point $y\in Y$, the  cokernel of the morphism $\pi^*: \bZ = K^*/R^* \to (M_{Y}^{gr}/\cO_Y^*)_y$ is torsion free.\\
 
 In (\cite{ikn}, Theorem 6.3 and Theorem 7.1) Illusie, Kato and Nakayama proved that under the above assumption the relative log de Rham cohomology $H^m( \oX,   \Omega ^*_{\oX /R} (log))$ is a free $R$-module, 
 the residue of the logarithmic Gauss-Manin connection on $H^m( \oX,   \Omega ^*_{\oX /R} (log))$  is nilpotent, 
 the Hodge spectral sequence, defined by the ``stupid'' filtration $\sigma_{\geq \bullet}$ on  $ \Omega ^*_{\oX /R} (log)$, degenerates at the $E_1$ term
 and the Hodge filtration 
 $$H^m( \oX,  \sigma_{\geq j}  \Omega ^*_{\oX /R} (log))  \mono H^m( \oX,   \Omega ^*_{\oX /R} (log)),$$
 splits ({\it i.e.} the associated graded $R$-module is free).
Moreover,  there is a canonical isomorphism\footnote{ The isomorphism (\ref{bhcomp})  depends on the choice of a uniformizer of $R$. Our choice is $t$.}
\begin{equation}\label{bhcomp}
 H^m(\tilde Y^{ log}, \bC) \simeq H^m( Y,   \Omega ^*_{Y/\bC} (log))
  \end{equation}
compatible with the base change $Y \rightsquigarrow Y_e$.
Set
$$F^j H^m( Y,   \Omega ^*_{Y/\bC} (log)) := H^m( Y,  \sigma_{\geq j}  \Omega ^*_{Y/\bC} (log)) \mono H^m( Y,   \Omega ^*_{Y/\bC} (log)).$$
As an immediate corollary of the Illusie-Kato-Nakayama results we get the following statement.
 \begin{lm}\label{indep2} \begin{enumerate}[(a)]  
 \item 
    Let $f: (\oX, M_{\oX})  \to (\oX', M_{\oX'})$ be a log morphism of  smooth  proper vertical fs log schemes over $(\mathrm{spec}\, R,  M_R)$ satisfying the condition (U).
  Assume that the induced morphism of generic fibers $f_K: \oX \otimes K \to \oX' \otimes K$ is an isomorphism. Then, for every non-negative integer $m$,  the morphism
   $$f^*: F^{\bullet} H^m( Y',   \Omega ^*_{Y'/\bC} (log)) \to  F^{\bullet} H^m( Y,   \Omega ^*_{Y/\bC} (log))$$
   is  a filtered isomorphism.
   \item
For a smooth  proper vertical fs log scheme $(\oX, M_{\oX})$  satisfying the condition (U) and a positive integer $e$ the canonical morphism 
  $$F^{\bullet} H^m( Y,   \Omega ^*_{Y/\bC} (log)) \to F^{\bullet} H^m( Y_e,   \Omega ^*_{Y_e/\bC} (log))$$ 
  is  a filtered isomorphism.
  \end{enumerate}
\end{lm}

 Assume, in addition, that $\oX$ is projective. Let  
 $$ W_{\bullet}= W_{\bullet} H^m(\tilde Y^{log}, \bQ)\subset H^m(\tilde Y^{log}, \bQ)$$
 be the monodromy filtration defined by the nilpotent endomorphism $N_Y= log\, T_Y$ of  $H^m(\tilde Y^{log}, \bQ)$:
$$N_Y W_i \subset W_{i-2},$$
$$Gr \, N^i_Y:   W_{i+m}^{\bQ}/W_{i+m-1}^{\bQ}  \iso  W_{m-i}^{\bQ}/W_{m-i-1}^{\bQ}.$$
 \begin{lm}\label{U} 
  For every smooth projective vertical fs log scheme $(\oX, M_{\oX})$  satisfying the condition (U)  the triple $(H^m(\tilde Y^{log}, \bZ), W_i H^m(\tilde Y^{log}, \bQ), F^j H^m( Y_{an},   \Omega ^*_{Y_{an}/\bC} (log)))$
  together with the isomorphism (\ref{bhcomp}) constitute a mixed Hodge structure.
 \end{lm}
\begin{proof}   By the semi-stable reduction theorem (\cite{kkms}, p. 198) we can find an integer $e$,  a projective strictly semi-stable scheme $\oX'$ over $R_e$ and log morphism
$$f: (\oX' , M_{\oX'}) \to (\oX_e, M_{\oX_e}),$$
where $M_{\oX'}$ is given by (\ref{can.log.str}) and $f$ is an isomorphism over the generic point of $R_e$.
The Lemmas \ref{indep} and \ref{indep2} reduce the proof to the case when  $\oX= \oX'$.  In  this case our assertion  is proven in (\cite{kawnam}, p. 405-406 and \cite{st3}, \S 5.6).
  \end{proof}

To construct the functor (\ref{fhs}) we define an auxiliary subcategory $SS_R$ of the category of schemes over $R$  whose objects are regular projective $R$-schemes $\oX$ such that the reduced special $Y_{red}$ is a strict normal crossing divisor on $\oX$.
 Let ${\cS} \subset Mor(SS_R)$ be the subset that consists of morphisms $f: \oX \to \oX'$  such that  $f_K: \oX \otimes K \simeq \oX' \otimes K$. 
\begin{lm}\label{localiz} 
The set $\cS$ is a left multiplicative system in $Mor(SS_R)$ (\cite{ks}, \S 7).
Moreover, the functor
 $$SS_{R}\to SmPr_K$$
 that  takes $\oX$ to  $ \oX \otimes  K $  exhibits  the category $SmPr_K$ as the localization of  $SS_{R}$ by $\cS$.
 \end{lm}
\begin{proof} The Lemma follows from the Hironaka theorem on resolution of singularieties  immediately. 
\end{proof}
Thus, by the universal property of the localization giving a functor from the category  $SmPr_K$ to another category is equivalent to giving a functor from $SS_{R}$ that takes every morphism in $\cS$ to an isomorphism. We define a functor
 $\Psi: SS_{R} \to  \widetilde {MHS}$ as follows.  Let $\oX$  be an object of $SS_{R} $, and let $M_{\oX}$ be the canonical log structure given by the formula (\ref{can.log.str}).  For sufficiently divisible integer $e$ the log scheme $(Y_e, M_{Y_e})$ satisfies the property (U). 
 We set
 \begin{align}\label{defpsi}
 \Psi( \oX)&=(H^m(\tilde Y^{log}, \bZ) \\
\nonumber &\simeq H^m(\tilde Y^{log}_e, \bZ), W_i H^m(\tilde Y^{log}_e, \bQ), F^j H^m( Y_e,   \Omega ^*_{Y_e/\bC} (log)), T_Y).
 \end{align}
 The right-hand side of (\ref{defpsi}) is independent of $e$ and is
   naturally promoted to a contravariant functor $\Psi: SS_{R} \to  \widetilde {MHS}$. By lemma  \ref{indep}  $\Psi$ takes every morphism in $\cS$ to an isomorphism.
   The functor (\ref{fhs}) is constructed. Let us check the required properties of (\ref{fhs}).

a) The only non-trivial statement is that $S$ preserves the Hodge filtration on $H^m( Y_{e},   \Omega ^*_{Y_e/\bC} (log))$.  Consider the  action of the group
$\bZ/e\bZ$ on the log scheme $(\oX_{e}, M_{\oX_{e}}) $ induced by the Galois action on $R_e$. The restriction of this action to $(Y_e, M_{\oX_{e}})$ yields  an action on $H^m( Y_{e},   \Omega ^*_{Y_e/\bC} (log))$.
One easily checks that the action of the generator $ 1\in \bZ/e\bZ$ on  $ H^m( Y_{e},   \Omega ^*_{Y_e/\bC} (log))$ equals  $S$. The compatibility with the Hodge filtration follows immediately.

b) This follows from (\ref{bc3}) and Lemma \ref{indep2}.

c)  Denote by $Y_{an}=Y(\bC)$ the analytic space associated with $Y$.  Let  $\tilde p: \tilde Y^{log} =Y^{log}\times_ {S^1} \bR^1 \to Y_{an}$  be the composition of the projection to the first factor and the canonical map  $p: Y^{log}\to Y_{an}$.
 The complex $R\tilde p_* \bZ \in D^b(Sh(Y_{an})) $ has a canonical automorphism $T_Y$ induced  by the automorphism of the space  $\tilde Y^{log}$. In (\cite{kawnam}, p. 405-406), Kawamata and Namikawa put 
 a weight filtration on the complex $R\tilde p_* \bQ$ and proved that this filtration yields the required spectral sequence with rational coefficients.  Thus, we just need to lift the Kawamata-Namikawa filtration to  $R\tilde p_* \bZ $.  The required canonical lifting is provided by the following result.
 
 \begin{lm}\label{perverse} ({\it cf. } \cite{sa}, Prop. 2.7) Assume that  $\oX$ is a strictly semi-stable scheme over $R$ of relative dimension $d$.
 \begin{enumerate}[(a)]  
 \item The complex $R\tilde p_* \bZ $ is a $(-d)$-shifted perverse sheaf on $Y$ ({\it i.e.}, $R\tilde p_* \bZ [d]$ is a perverse sheaf).  Moreover,  the canonical filtration $\tau_{\leq i} R\tilde p_* \bZ  $ is a filtration by $(-d)$-shifted perverse subsheaves and it coincides with the filtration
 on $R\tilde p_* \bZ $ by kernels of $(T-1)^{i+1}$ (computed in the abelian category of $(-d)$-shifted perverse sheaves):
 $$\tau_{\leq i} R\tilde p_* \bZ = \Ker\left((T-1)^{i+1}:   R\tilde p_* \bZ \to R\tilde p_* \bZ\right).$$
  In particular,    $(T-1)^{d+1}$ is $0$ on $ R\tilde p_* \bZ $.
\item Let $0\subset  W_{-d} R\tilde p_* \bZ \subset \cdots  W_i  R\tilde p_* \bZ \subset  W_{d} R\tilde p_* \bZ = R\tilde p_* \bZ  $ be the monodromy filtration on  $R\tilde p_* \bZ  $ viewed as an object of the abelian category of $(-d)$-shifted perverse sheaves equipped with the nilpotent endomorphism $T-1$.
Then, for every integer $r$, we have an isomorphism 
$$Gr_r^W R\tilde p_* \bZ \simeq  \bigoplus_{\stackrel{i-j=r}{ i,j\geq 0}}  \wedge^{i+j+1} (M^{gr}_{Y_{an}}/\cO_{Y_{an}}^*)\simeq  \bigoplus_{i-j=r} a_{i+j  *} \bZ[-i-j],$$
where $a_q$ denotes the embedding $Y^{(q)}\mono Y$. The first isomorphism is canonical, the second one depends on the order of the set of irreducible components of $Y$.
\item The Verdier dual complex $\cD_Y(R\tilde p_* \bZ )$ is quasi-isomorphic to $R\tilde p_* \bZ [2\dim \, X]$.
\end{enumerate}
  \end{lm}
  \begin{proof} For the first statement it suffices to prove that, for every prime number $\ell$,  the complex $\bZ_{\ell} \otimes  R\tilde p_* \bZ=R\tilde p_* \bZ_{\ell}$ has the corresponding properties.  According to the comparison results of Kato and Nakayama (see, e.g. \cite{il2}, Th. 5.9,  Cor. 8.4)  the complex $R\tilde p_* \bZ_{\ell}$ is quasi-isomorphic 
 to the complex of nearby cycles $R\Psi \bZ_{\ell}$ computed using the \'etale  topology.  The results for  $R\Psi \bZ_{\ell}$ are proven in  (\cite{sa}, Lemma 2.5 and Cor. 2.6). The proof of the second statement is parallel to the proof of the analogous result for $R\Psi \bZ_{\ell}$  ( \cite{sa}, Prop. 2.7).
 For the last statement of the Lemma observe that $R\tilde p_* \bZ $ is quasi-isomorphic to $R p_{1*} \bZ$, where $p_1:  Y^{log}_1\to Y_{an}$ is the restriction of the map $Y^{log}\to   Y_{an} \times S^1$ to the fiber over $Y_{an} \times \{1\}$. As the map $p_1$ is proper, we have
 $$\cD_Y( R p_{1*} \bZ)\simeq Rp_{1*} \cD_{Y^{log}_1} \bZ.$$ 
Finally, a simple local computation shows that $\cD_{Y^{log}_1} \bZ \simeq \bZ[2 \dim \, X]$.
\end{proof}  
  d) The last assertion of Theorem \ref{lhs} is proven in  (\cite{ikn} Theorem 8.3, \cite{st2}, Appendix). The proof of Theorem \ref{lhs} is now completed.

\end{proof}  
  %\begin{rem}
 % The mixed Hodge structure $H^m(\lim X)$ depends on the choice of coordinate $t$  but the pure $\bQ$-Hodge structures $Gr_W^i H^m(\lim X) \otimes \bQ$ do not.
  %\end{rem}
  \begin{rem} We expect that the functor (\ref{fhs}) extends to the category of smooth quasi-compact rigid analytic varieties over $\bC((t))$ ({\it cf.} \cite{a2}).
 \end{rem}

We finish this subsection by recalling a variant of the Picard-Lefschetz formula for semi-stable degenerations.
Let $\oX$ be a projective strictly semi-stable scheme over $R$ of relative dimension $d$, and let $Y$ be its special fiber.  The simplicial complex that computes the homology of the Clemens polytope $Cl(Y)$ coincides with the complex
$$E_1^{-d, 2d}(\oX)(d) \to E_1^{-d+1, 2d}(\oX)(d)\to \cdots \to E_1^{0,2d}(\oX)(d),$$
where $E_r^{pq}(\oX)$ is the weight spectral sequence from Theorem \ref{lhs}.
From this we get a canonical morphism
\begin{equation}\label{clw0}
 H^m(\lim X, \bZ) \to  E_2^{-d+m, 2d}(\oX)(d) \simeq H_m(Cl(Y))  
  \end{equation}
  As the weight spectral sequence degenerates rationally in $E_2$ terms (\ref{clw0}) yields an isomorphism 
  \begin{equation}\label{clw0r}
  H_m(Cl(Y))\otimes \bQ  \stackrel{\gamma}{\simeq} Gr_{2d}^{W^\bQ}H^m(\lim X)(d).
  \end{equation}
 We apply this to $m=d$. If $$<, >:  W^\bQ_0 H^d(\lim X)   \otimes  Gr_{2d}^{W^\bQ}H^d(\lim X)(d)  \to \bQ$$ denotes the pairing induced by Poincar\'e duality (Lemma \ref{perverse}, c))
 then, for every 
 $$x=\sum_{v\in \pi_0(Y^{(d)})} a_v v, \quad  y=\sum_{v\in \pi_0(Y^{(d)})} b_v v \in H_d(Cl(Y))\otimes \bQ, $$  we have
   \begin{equation}\label{picle}
   (-1)^{\frac{d(d-1)}{2}}  <Gr N^d \gamma(x),  \gamma(y)>=  \sum_va_v b_v.
   \end{equation}
 This follows from compatibility of the weight spectral sequence with Poincar\'e duality and the monodromy action(\cite{sa}, Cor. 2.6 and Prop. 2.15).   
   % We expect that $H^m(\lim X)$ carries a canonical {\it integral} weight filtration. That is, for every morphism  $(f, f'):   (e_1,  \oX'  \rar{h_X}     \oX) \to ( e_2,  \oZ'  \rar{h_Z}     \oZ)$
 %in $\cS$, the induced morphism  
 %$$H^m(\tilde Z^{\prime log}, \bZ) \to H^m(\tilde X^{\prime log}, \bZ)$$
 %is strictly compatible with the integral weight filtration defined in  (\ref{integralweight}). However, we can not prove this.
% \subsection{Motivic integral of analytic Calabi-Yau varieties over a non-archimedean field.}\label{m.i.c.y.v} 
%Although in this paper we mostly interested in motivic integrals of algebraic varieties in some proofs  (see, {\it e.g.} \S \ref{k.k3})  we have to consider a broader class of space.  In this subsection we review the construction of motivic integral of 
%analytic (either in the sense Berkovich (\cite{ber1}) or in the sense of rigid Geometry)  Calabi-Yau varieties over a non-archimedean field.
%\subsection{Motivic integral of Abelian varieties over $\bC((t))$.}\label{m.i.a.v.}
\subsection{Motivic Serre Invariant.}\label{m.s.i.}  
  Let $R$ be a complete discrete valuation ring with perfect residue field $k$ and fraction field $K$.
  The motivic Serre invariant of  a smooth proper variety $X$ over $K$ is the class of the special fiber $V^0$ of a weak N\'eron model $\cV$ of $X$  
 in the quotient ring
  $$K_0(Var_k)_{loc}\to K_0(Var_k)_{loc}/ (\bZ(1)-\bZ).$$
   It is shown in (\cite{ls}, Theorem 4.5.1) that the motivic Serre invariant $S(X)$ is well defined  {\it i.e.},  independent of the choice of $\cV$.
   If $X$ is a Calabi-Yau variety $S(X)$ equals the image of the motivic integral $\int_{X}$ in the  quotient ring.
   
   Let $K= \bC((t))$. In the following Proposition, which is a refinement of  A'Campo's formula for the Euler characteristic of the motivic integral\footnote{Related results were obtained by Nicaise (\cite{ni}).}, we denote by  $S^H(X)$ the image of $S(X)$  under the ring homomorphim
    \begin{equation}\label{virth}
    K_0(Var_{\bC})_{loc}/ (\bZ(1)-\bZ) \to K_0(MHS)/ (\bZ(1)-\bZ)
    \end{equation}
    that takes the class of a variety $Z$ to the virtual mixed Hodge structure   $\sum (-1)^{i}[H_c^i(Z,\bZ)]$. 
         
\begin{pr}\label{modQ(1)-Q}
Let $X$ be a smooth projective variety over $\bC((t))$. Assume that $X$ has  a projective strictly semi-stable model $\oX$ over $\bC[[t]]$.  Then  $S^H(X)$ 
 is equal to the class of $\sum(-1)^i [H^i(\lim X)]$. 
\end{pr}
\begin{proof} We start with the following general (and well known)  observation. 
\begin{lm}\label{exampneron} Let $R$ be a  complete discrete valuation ring with perfect residue field $k$ and fraction field $K$, and let $\oX$ be a proper flat scheme over $R$. Assume that $\oX$ is regular and that the generic fiber
$X= \oX\otimes _R K$ is smooth over $K$.  Then the smooth locus $\oX_{sm}$ of the morphism $\oX\to \mathrm{spec}\, R$ is a weak N\'eron model of $X$.
\end{lm}
 \begin{proof} Since $X$ is smooth we have that  $\oX_{sm}\otimes _R K = X$. Let $R'\supset R$ be a finite unramified extension with fraction field $K'$. We need to show that every
morphism $x: \mathrm{spec} \, K' \to X$ extends to  an $R$-morphism $\overline x: \mathrm{spec} \, R' \to \oX_{sm}$. As $\oX$ is proper over $R$, $x$ extends to an $R$-morphism $\overline x: \mathrm{spec} \, R' \to \oX$. We claim that $\overline x$ 
 takes the closed point of $\mathrm{spec}\, R'$ to a smooth point, $y$, of the special fiber $Y=\oX\otimes _R k$. Since $k$ is perfect, it suffices to check that $y$ is a regular point of $Y$ (\cite{sga1}, {\rm II}, Cor. 5.3). 
 Indeed, let $\cO_{\oX, y}$ (resp.  $\cO_{Y, y}$) be the local ring of $\oX$ (resp. $Y$) at $y$ and let $m_{\oX, y}\subset \cO_{\oX, y}$ (resp. $m_{Y, y}\subset \cO_{Y, y}$) be the maximal ideal. We have a surjective morphism
 \begin{equation}\label{neproptr}
 m_{\oX, y}/m_{\oX, y}^2 \epi m_{Y, y}/m_{Y, y}^2
 \end{equation}
 of finite-dimensional vector spaces over $\cO_{\oX, y} /m_{\oX, y}$. Let us show that the image in  $m_{\oX, y}/m_{\oX, y}^2$ of a uniformizer $t\in R$  is not equal to $0$. Indeed, we have
 a morphsim $ \cO_{\oX, y} \rar{\overline x^*}  R'$ induced by $\overline x$ such that the composition $R\to \cO_{\oX, y} \rar{\overline x^*}  R'$ is the identity morphism. 
 Since $K'$ is unramified over $K$, $t$ is also a uniformizer for $R'$. Therefore, $t$ does not belong to $m_{\oX, y}^2$. We proved that the image of $t$ in $m_{\oX, y}/m_{\oX, y}^2$ is not $0$. On the other
 hand, its image in $m_{Y, y}/m_{Y, y}^2$ is $0$. Hence, morphism (\ref{neproptr}) is not injective and, therefore,
 $$\dim  m_{\oX, y}/m_{\oX, y}^2  > \dim m_{Y, y}/m_{Y, y}^2.$$
 On the other hand, since $\oX$ is regular, we have that
   $\dim  m_{\oX, y}/m_{\oX, y}^2$ equals the Krull dimension of  $\cO_{\oX, y}$. Thus,  $\dim m_{Y, y}/m_{Y, y}^2\leq \dim  \cO_{\oX, y} -1 = \dim  \cO_{Y, y}$. Hence, $Y$ is regular and, therefore, smooth at point $y$. It follows that the map  $\overline x: \mathrm{spec} \, R' \to \oX$ factors through $\oX_{sm}\subset \oX$.
  \end{proof}
We now come back to the proof of Proposition \ref{modQ(1)-Q}. According to the above lemma the smooth locus $\cV$ of $\oX$ is a weak N\'eron model of $X$. 
Using notation of (\ref{sstm}) and the inclusion-exclusion formula we find  
$$[V^0]=  \sum_{j=0}^{\dim\, X} \left((-1)^j(j+1)\left[Y^{(j)}\right] \right).$$
 On the other hand, by part c) of the Theorem \ref{lhs} the class $\sum(-1)^i [H^i(\lim X)]$ is equal to the image under (\ref{virth}) of the class
 \begin{align*}
[\lim X]=\sum_{j=0}^{\dim\, X}\left((-1)^j  \left[Y^{(j)}\right] \sum_{a=0}^{j}\bZ(-a)\right).
\end{align*} 
 Comparing the two formulas we complete the proof  of Proposition \ref{modQ(1)-Q}.
 \end{proof} 
Let $\chi:K_0(Var_{\bC})\to \bZ$ be the ring homomorphism defined by
$$\chi\left([Z]\right)=\sum (-1)^i\dim H^i_c(Z,\bC).$$
Notice that since $\chi(\bZ(1)-\bZ)=0$, $\chi$ factors uniquely through $ K_0(Var_{\bC})_{loc}/ (\bZ(1)-\bZ)$. We have the following corollary of Proposition \ref{modQ(1)-Q}.
\begin{cor}[cf. A'Campo (\cite{ac})]\label{acampo}
Let $X$ be a smooth projective variety over $K=\bC((t))$. Assume that $X$ has a projective strictly semi-stable model $\overline{X}$ over $\bC[[t]]$. Then
$$\chi\left(S(X)\right)=\sum(-1)^i\dim H^i(\lim X, \bC).$$
\end{cor}

 In the rest of this subsection, we explain an analogue of the above Proposition for the finite residue field case. Let $K$ be a local field with residue field $k=\bF_q$, and let
 $$K_0(Var_{\bF_q})_{loc}/(\bZ(1)-\bZ) \to \bZ/(q-1), $$
 be the homomorphism  induced by (\ref{p.real}).  The image of $S(X)$ in  $\bZ/(q-1)$ is the classical Serre invariant which we denote by $S^q(X)$. 
  \begin{pr}\label{modq-1}
Let $X$ be a smooth proper variety over  $K$. Assume that $X$ has a proper strictly semi-stable model over the ring of integers $R$. 
Then the Serre invariant of $X$ is given by the formula
\begin{equation}\label{m.q}
 \sum _j Tr(F^{-1}, H^j(X_{\oK}, \bQ _{\ell}) 
\end{equation}
where $F\in \mathrm{Gal} (\oK/K)$ is a lifting of the Frobenius automorphism  $Fr\in \mathrm{Gal}(\overline k/k)$ and $\ell$ is a prime number different from the characteristic of $k$.  
 \end{pr} 
  \begin{proof}
  This can be proved as its Hodge analogue above using the $\ell$-adic weight spectral sequence. We give a different proof.
   Let $\oX$ be  a strictly semi-stable model of $X$. Then the Serre invariant of $X$ equals $|Y_{sm}(k)|$ modulo $(q-1)$. 
On the other hand, if $\Psi(\bQ _{\ell})$ is the complex of nearby cycles (viewed as a complex of $\ell$-adic sheaves on $Y$), 
by the Grothendieck-Lefschetz formula we have
 \begin{equation}\label{grltrfor}
  \sum _j (-1)^jTr(F^{-1}, H^j (X_{\oK}, \bQ _{\ell}))= \sum _j (-1)^jTr(F^{-1}, H^j (Y_{\overline k}, \Psi(\bQ _{\ell})))=$$
$$\sum_{y\in Y(k)} \sum _i (-1)^iTr(F^{-1}, \cH^i ( \Psi(\bQ _{\ell}))_y).
\end{equation}
If $y\in Y_{sm}(k)$, the corresponding internal sum equals $1$. 
If  $y\in Y_{sing}(k)$ then $ \cH^i ( \Psi(\bQ _{\ell}))_y\simeq \bigwedge^i T (-i)$, where $T$ is a vector space with the trivial
action of $\mathrm{Gal}(\oK/K)$ (\cite{sga7}, Expos\'e {\rm I}, Th. 3.3). Thus, for $y\in Y_{sing}(k)$, we have
$$  \sum _i (-1)^iTr(F^{-1}, \cH^i ( \Psi(\bQ _{\ell}))_y)\equiv  \sum _i (-1)^i \dim\, \bigwedge^i T \equiv   0 \mod(q-1).$$
It follows that the right-hand side of (\ref{grltrfor}) is equal to $|Y_{sm}(k)|$ modulo $(q-1)$ which is the Serre invariant of $X$.
 \end{proof}
 
  \section{Motivic integral of K3  surfaces over $\bC((t))$ .}\label{p.m.t}
  In this section we will prove Theorem \ref{main} stated in the introduction. Without loss of generality we may assume that the ramification index $e$ is equal to 1. 
  Indeed, by Theorem \ref{lhs} part (b), the formulas (\ref{eq1m}) and  (\ref{eq2m}) for the pair $(X/K, e)$ are equivalent to those for the pair  $(X_{K_e}/K_e, 1)$. If $X$ admits a strictly semi-stable model over $R$ then  $X_{K_e}$ admits  a strictly semi-stable model over $R_e$ (\cite{sa}, Lemma 1.11).
  We will write $r_i$ for $r_i(X,K)$.
 \subsection{Approximation of varieties over the formal disk.}\label{appr} We will need the following version of Artin's Approximation Theorem. 
\begin{pr}\label{artinappr} Let $k$ be a field of characteristic $0$, and let $\oX$ be a projective strictly semi-stable scheme over $R= k[[t]]$. For every positive integer $n$ there exist
\begin{enumerate}[(1)]
\item  a smooth curve $\oC$ over $k$ with a point $a\in \oC(k)$,
\item an \'etale morphism $h: \oC \to \bA^1_k= \mathrm{spec}\, k[t]$ that carries $a$ to $0$,
\item a flat projective scheme $\overline \cX$ over $\oC$,
\item an isomorphism of schemes over $R_n= \mathrm{spec}\, k[t]/t^{n+1}$:
$$\oX \times_{\mathrm{spec}\, R } \mathrm{spec} \, R_n \simeq  \overline \cX \times_\oC \mathrm{spec} \, R_n.$$
\end{enumerate}
Here $\mathrm{spec} \, R_n$ is viewed as a scheme over  $\oC$ via the unique morphism $\tilde i_n: \mathrm{spec} \, R_n \to \oC$ that carries the point $0$ to $a$ and makes the following diagram commutative
$$
\xymatrix{
& \oC\ar[d]_{h}    \\ 
\mathrm{spec} \,  R_n \ar@{^{(}->}[r]_{i_n} \ar[ur]_{\tilde i_n}& \mathrm{spec} \, k[t]   \\
   }
$$

If $\oC, a, h,  \overline \cX $ are as above, the scheme  $\overline \cX $ is  regular in an open neighborhood of its special fiber $Y'$ and $Y'$ is a reduced divisor on $\overline \cX $
with strict normal crossings. In addition, if $X$ is a $d$-dimensional Calabi-Yau variety
the collection  $\oC, a, h,  \overline \cX $ can be chosen so that the line bundle $\Omega^d_{\cX/C}$ is trivial and
\begin{equation}\label{eqaprmotintegral}
\int_{X} =\int_{\cX \times_C \mathrm{spec}\, K'}.
\end{equation}
 Here we set $C=\oC-a$, $\cX=  \overline \cX \times_\oC C$, and  $K'$ denotes the fraction field of the completed local ring  $R'=\hat \cO_{\oC,a}$. 
\end{pr} 
\begin{proof}
Choose an embedding $\oX\mono \bP_{R}^n$ and let $\nu:\mathrm{spec}\ R \to \mathrm{Hilb}(\bP^n_{R})$ be the corresponding morphism to the Hilbert scheme. Using Artin's Formal Approximation Theorem (see {\it e.g.}, \cite{blr}, \S 3.6) on the morphism $\nu$ we obtain (1)-(4). Next, 
 we claim that the scheme $\overline{X}' =\overline{\cX}\times _\oC \mathrm{spec}\, R'$ is regular. As  $\overline{X}'$ is proper over $R'$ and the set of its regular points is open (\cite{ega} {\rm IV}, 6.12.5) 
 it suffices to show that the local ring of any point of the special fiber $Y'$ is regular which in turn follows from property (4) and the regularity of $\oX$.  Moreover, $Y'$ being isomorphic to the special
 fiber of a strictly semi-stable scheme $\oX$ is a strict normal crossing divisor on $\oX'$ and on $\overline{\cX}$.  
 Note that under our assumption that $\mathrm{char}\,k =0$ this implies strict semi-stability of $\oX'$.
% Now let $b$ be any point of the special fiber and let $m\subset \cO_{\overline{\cX},b'}$ be the maximal ideal of $\cO_{\overline{\cX},b'}$, then $m/m^2$ has the same dimension as $\cX$ by part (4) and the fact that $\dim \cO_{Y,b'}=\dim \cO_{\cX,b'}-1.$ 
%Since the total space is regular, $\mathrm{char}\ k=0$, $\dim \oX=\dim \overline{\cX}$, $\oX$ was assumed to be strictly semi-stable it follows that $\overline{\cX}$ satisfies the properties in section \ref{c.p.} and hence %is strictly semi-stable over $R'$. Note that $\mathrm{char}\ k=0$ is needed for the regularity of the total space to imply the smoothness of the generic fiber.

 %It remains to prove the last assertion of the proposition. 
 Suppose that $X$ is a Calabi-Yau variety. Then the divisor of any nonzero relative log form $\omega \in H^0(\oX,\Omega^d_{\oX/R}(log))$ 
 is supported on the special fiber $Y$ of $\oX$. Write $\mathrm{div}(\omega)=\sum_i n_i[V_i]$, 
 where $V_i$ are the irreducible components of $Y$.  Assume that  the quadruple $\oC, a, h,  \overline \cX $ 
 satisfies properties (1)-(4) with $n\ge \sum_{i}{n_i}$. To prove the last assertion of the proposition, formula  (\ref{eqaprmotintegral}), we will show that there exists a section $\omega' \in H^0(\oX',\Omega^d_{\oX'/R'}(log))$ whose divisor is supported on the special fiber $Y'$ of $\oX'$ and 
 such that via the isomorphism $Y\simeq Y'$ from (4) 
 \begin{equation}\label{eqofdiv}
 \mathrm{div}(\omega)=\mathrm{div}(\omega').
 \end{equation}
 Indeed, by Lemma 4.1 from  (\cite{kawnam}), for every proper strictly semi-stable scheme $\oX$ over $R$ the $R$-module $H^0(\oX,\Omega^d_{\oX/R}(log))$ is free and, in addition,  we have
 $$H^0(\oX,\Omega^d_{\oX/R}(log))\otimes_R R_n \iso H^0(\oX\otimes R_n,\Omega^d_{\oX\otimes R_n /R_n}(log)).$$ 
 Applying this result to $\oX$ and $\oX'$ we find that  $H^0(\oX,\Omega^d_{\oX/R}(log))$ and $H^0(\oX',\Omega^d_{\oX'/R'}(log))$ are free modules of rank 1 over $R$ and $R'$ respectively
  and that
 (4) induces an isomorphism 
 $$\theta:H^0(\oX,\Omega^d_{\oX/R}(log))\otimes _{R} R_n \iso H^0(\oX',\Omega^d_{\oX'/R'}(log))\otimes _{R} R_n.$$
 (The $R$-action on  $H^0(\oX',\Omega^d_{\oX'/R'}(log))$ comes via the isomorphism $R\iso R' $ induced by $h$.) We claim that a section $\omega' \in H^0(\oX',\Omega^d_{\oX'/R'}(log))$ such that 
 $\theta(\omega \otimes 1)=\omega'\otimes 1$ does the job. Our claim is local: it suffices to show that, for a closed point $b\in \oX$ and local regular functions $f,g\in \cO_{\oX, b} $  such that
 $\mathrm{div}(f)$ is supported on $Y$, $\sum_i \mathrm{ord}_{V_i} f \leq n$, and $f-g \in( t^{n+1})$, one has
  $\mathrm{div}(f)=\mathrm{div}(g)$.  Let $x_i$ be a system of local parameters at $b$ such that $t=x_1\cdots x_m$.  Then, locally around $b$, we have  $f=x_1^{n_1}\cdots x_m^{n_m}u$, where $u$ is invertible and $\sum_i n_i\leq n$.
  If $n_1>0$, $g\in f+ (t^{n+1})$ is divisible by $x_1$ and $\frac{f}{x_1}-\frac{g}{x_1}\in (t^n)$. Arguing by induction we see that $g$ is divisible by $x_1^{n_1}\cdots x_m^{n_m}$
  and
  $$\frac{f}{x_1^{n_1}\cdots x_m^{n_m}}-\frac{g}{x_1^{n_1}\cdots x_m^{n_m}}\in (t).$$
  In particular, $g=x_1^{n_1}\cdots x_m^{n_m}u'$ for some invertible $u'$.
  
  To complete the proof of the proposition let us explain how (\ref{eqofdiv}) implies (\ref{eqaprmotintegral}). Suppose that the pair $ \overline \cX $, $\omega' \in H^0(\oX',\Omega^d_{\oX'/R'}(log))$ is chosen such that
   that equation (\ref{eqofdiv})  holds. Then, in particular, $\omega'$ restricts to a non-vanishing differential form on $X'$. Thus, $X'$ is a Calabi-Yau variety. Secondly, by  Lemma \ref{exampneron} the schemes
   $\oX_{sm}$ and $\oX'_{sm}$ are weak N\'eron models of $X$ and $X'$ respectively. Moreover, by property (4) and (\ref{eqofdiv}) there exists an isomorphism between the special fibers of  $\oX_{sm}$ and $\oX'_{sm}$
   that carries $\mathrm{div}(\omega)$ to $\mathrm{div}(\omega')$. Using (\ref{n.m.i}) formula (\ref{eqaprmotintegral}) follows. 
    \end{proof}
 \subsection{Kulikov model.}\label{ss.kulikovmodel}
 It is enough to prove Theorem \ref{main} in the case where $X$ is the restriction of a strictly semi-stable family over a complex curve. Indeed, apply Proposition \ref{artinappr} to a strictly semi-stable model $\oX$ of $X$. As the limit mixed Hodge structure of a strictly semi-stable scheme depends only
 on its special fiber together with its log structure which, in turn, is determined by its first infinitesimal neighborhood  $\oX\otimes _R R/t^2$, the formulas (\ref{eq1m}),  (\ref{eq2m}) for $X$ are equivalent to those for $\cX \times_C \mathrm{spec}\, K'$. 
 
  Let $X$ be a K3 surface over $K$, which is the restriction of a strictly semi-stable family over a complex curve.  In (\cite{ku}, Theorem 2),  Kulikov demonstrated that $X$ has a projective strictly semi-stable model $\oX$ over $R$ such that 
  the log canonical bundle $\Omega^2_{\oX/R}(log)$ is trivial and
  the special fiber $Y$ is of one of the following types (depending on the number $s$  defined in Theorem \ref{main})
 \begin{enumerate}[(I)]
 \item ($s=1$)  $Y$ is a smooth K3 surface
 \item ($s=2$)   $Y$ is a chain of smooth surfaces $V_0,\ldots,V_m$ ruled by elliptic curves, with smooth rational surfaces on either end and each double curve $V_{i}\cap V_{i+1}$ is a smooth elliptic curve. 
 \item  ($s=3$) $Y$ is a union of smooth rational surfaces whose pairwise intersections are smooth rational curves and the Clemens polytope of $Y$ is a triangulation of $S^2.$
 \end{enumerate}
 In addition, for $s=2$, Friedman showed in (\cite{fr}, Theorem 2.2) that a Kulikov model can be chosen so that all the ruled elliptic surfaces in $Y$ are minimal {\it i.e.}, $\bP^1$-fibrations over an elliptic curve. We shall call such model {\it special}.     

If $\oX$ is a Kulikov model, we have
 \begin{equation}\label{minkul}
 \int_{X}=[Y_{sm}].
  \end{equation}
 Indeed, by Lemma \ref{exampneron} the smooth locus $\oX_{sm}$ of $\oX$ is a weak N\'eron model of $X$. Moreover, since the log canonical bundle $\Omega^2_{\oX/R}(log)$ is trivial, the bundle
 $\Omega^2_{\oX_{sm}/R}$ (which is isomorphic to the restriction of  $\Omega^2_{\oX/R}(log)$ to $\oX_{sm}$) is also trivial.  If $\omega\in \Gamma(X_{sm}, \Omega^2_{\oX_{sm}/R})$ is a trivializing section,
 the numbers $m_i$ appearing in formula (\ref{def.mult}) are all equal to $0$.  Thus, by formula (\ref{n.m.i}) the motivic integral $\int_{X}$ is equal to the sum of classes of the irreducible components of $Y_{sm}$. Since $Y_{sm}$ is smooth
 its irreducible components are pairwise disjoint and, hence, the sum of its classes is equal to $[Y_{sm}]$.
  \subsection{Type II Degeneration}
Suppose that $\oX$ is a type II special Kulikov model.  Let $V_0,\ldots,V_m$ be the irreducible components of $Y$ such that $V_0$ and $V_m$ are rational surfaces, and let  $C_i=V_{i}\cap V_{i+1}$ be the double curves.
 \begin{lm}\label{uh} 
 \begin{enumerate}[(1)]
\item  Let $E_1,\ldots,E_{m-1}$ be ruling elliptic curves for $V_1,\ldots,V_{m-1}$.  Then $C_{i}\cong E_i\cong E_j\cong C_j$ for all $i$ and $j$. 
\item At least one of the rational components, $V_0$ or $V_m$, is not minimal.
\end{enumerate}
 \end{lm} 
\begin{proof}
(1). We will first prove that $E_1\cong E_2$. Let $C_1$ and $C_2$ be elliptic curves given by the intersection $V_1\cap V_2$ and $V_2\cap V_3$ respectively. We have the following diagram
$$
\xymatrix{
C_1 \ar@{^{(}->}[r]\ar[dr]_{f_1}& V_2\ar[d]_{h} & \ar@{_{(}->}[l]\ar[dl]^{f_2} C_2\\
             & E_2 &  \\ 
}
$$
Notice that the maps $f_1$ and $f_2$ cannot be constant since this would imply the existence of injections of $C_1$ and $C_2$ into rational curves. Thus $f_1$ and $f_2$ must be finite. 
The triviality of the log canonical bundle $\Omega^2_{Y}(log)$ implies that for the canonical class $K_{V_2}$ we have $K_{V_2}=-[C_1]-[C_2]$. On the other hand,  the restriction of $K_{V_2}$ to a smooth fiber, $h^{-1}(a)$, of the map $h: V_2\to E_2$ is isomorphic to $K_{h^{-1}(a)}$.  As $h^{-1}(a)$ is a smooth rational curve, we have that  $\deg(K_{V_2}|_{h^{-1}(a)})=-2$ which implies the degree of the divisor $-[C_1]-[C_2]$ intersected with the fiber $h^{-1}(a)$ is $-2$. Hence the images of $C_1$ and $C_2$ in $V_2$ have only one intersection point with a generic fiber which implies $f_1$ and $f_2$ are one-to-one and $C_1\cong E_2\cong C_2.$ We then apply the same method of proof to show that $C_2\cong E_3\cong C_3$ and so on.

(2). We claim that for a minimal ruled elliptic surface $V_i$ and two disjoint sections $C_{i-1}, C_i \subset V_i $, we have
 $$\left( [C_{i-1}]\right)^2_{V_i} = - \left([C_{i}]\right)^2_{V_i}.$$
 Indeed, the N\'eron-Severi group of $V_i$ is generated by the class $[C_i]$ of $C_i$ and the class $[\bP^1]$ of a smooth fiber of the map $V_i\to E_i$. If $[C_{i-1}]= [C_i] +c [\bP^1]$, we have
$$0= \left(c[\bP^1]\right)^2_{V_i}= \left( [C_{i-1}] - [C_i]\right)^2_{V_i}= \left( [C_{i-1}]\right)^2_{V_i} + \left([C_{i}]\right)^2_{V_i}.$$
On the other hand, since $Y$ is the special fiber of a semi-stable degeneration, we have for every $i$
 $$ \left([C_{i}]\right)^2_{V_i} = -\left([C_{i}]\right)^2_{V_{i+1}}.$$
 Combining the two formulas we see that  $\left([C_{0}]\right)^2_{V_0}=  -\left([C_{m-1}]\right)^2_{V_m}$. In particular,  at least for one of the rational components, say $V_0$, the self-intersection of the double curve lying on it is non-positive.  Thus, $(K_{V_0})^2_{V_0}= \left(-[C_{0}]\right)^2_{V_0}\leq 0$. Using Noether's formula (\cite{bea} I.14) it follows that
  $V_0$ is not minimal.
 \end{proof}
Let $E$ be an elliptic curve such that $E\cong C_i$ for all $i$. Then we get from (\ref{minkul}) 
\begin{align*}
\int_{X}&=\sum_{i=0}^{m}[V_i]-2\sum_{i=0}^{m-1}[C_{i}]=\sum_{i=0}^{m}[V_i]-2m[E].
\end{align*}
Since $V_0$ and $V_m$ are both rational surfaces we have $[V_0]=\bZ+a_0\bZ(-1)+\bZ(-2)$ and $[V_m]=\bZ+a_m\bZ(-1)+\bZ(-2)$. Each $V_i$ for $1\le i\le m-1$ is birationally equivalent to $\bP^1\times E$. Thus, by (\cite{bea} II.11), $[V_i]=[E\times \bP^1]+a_i\bZ(-1)$  for $1\le i\le m-1$.
 Letting $a=\sum_{i=0}^{m}a_i$ we have 
\begin{align*}
\int_{X}&=2\bZ + a\bZ(-1)+(m-1)[E]\cdot[\bP^1]+2\bZ(-2)-2m[E]\\
&=2\bZ + a\bZ(-1)+(m-1)[E]+(m-1)[E](-1)+2\bZ(-2)-2m[E]\\
&=2\bZ + a\bZ(-1)-(m+1)[E]+(m-1)[E](-1)+2\bZ(-2).
\end{align*} 
Using Corollary \ref{acampo} and the fact that the Euler characteristic of a K3 surface is 24 it follows that $a=20$. Thus we have the formula 
$$\int_{X}=2\bZ -(m+1)[E]+20\bZ(-1)+(m-1)[E](-1)+2\bZ(-2).$$
Now we want to express the number of double curves $m$ and the class of the elliptic curve $[E]$ in terms of the limit mixed Hodge structure $H^2(\lim X)$. 
%We claim that $m$ is equal to the order of the cokernel of the map $Gr\ N:W^{\bZ}_3/W^{\bZ}_1\to W^{\bZ}_1/W^{\bZ}_0$ and the class $[E]$ is defined by the rank 2 Hodge structure on $W_1^{\bZ}=W_1^{\bQ}\cap H^2(\lim X,\bZ)$. 
First, we show that the {\it integral}  weight spectral sequence $E^{pq}_r$ from Theorem \ref{lhs} degenerates at the second term. Since it degenerates rationally it will suffice to show that the $E_2$ terms are torsion free. The nontrivial portion of the first term of the spectral sequence is 
\begin{equation}\label{spectral}
\begin{aligned}
\xymatrix{
 \oplus_{i=0}^{m-1}H^2(C_i)(-1)\ar[r]^{\delta_4}&\oplus_{i=0}^{m}H^4(V_i) & \\
\oplus_{i=0}^{m-1}H^1(C_i)(-1) \ar[r]^{\delta_3}&\oplus_{i=1}^{m-1}H^3(V_i) & \\
\oplus_{i=0}^{m-1}H^0(C_i)(-1)\ar[r]^{\delta_2} &\oplus_{i=0}^{m}H^2(V_i)\ar[r]^{\delta_2'} &\oplus_{i=0}^{m-1}H^2(C_i) \\
 & \oplus_{i=1}^{m-1}H^1(V_i)\ar[r]^{\delta_1}& \oplus_{i=0}^{m-1}H^1(C_i)\\
 & \oplus_{i=0}^{m}H^0(V_i)\ar[r]^{\delta_0}&\oplus_{i=0}^{m-1}H^0(C_i) 
}
\end{aligned}
\end{equation}
The first and the last complexes compute (co)homology of the Clemens polytope of $Y$ and, hence, are quasi-isomorphic to $\bZ$.
Consider the middle complex. The map $\delta_2$ is injective since $\delta_2\otimes \bQ$ is. Let us prove that $\delta_2'$ is surjective.
For every $(u_0,\ldots,u_{m})\in \oplus_{i=0}^{m}H^2(V_i)$, we have
$$\delta_2'(u_0,\ldots,u_{m})=\left((u_0)_{|C_0} -(u_1)_{|C_0},\ldots, (u_{m-1})_{|C_{m-1}}- (u_m)_{|C_{m-1}}\right).$$
For every $ 1\leq i\leq m-1 $ the restriction morphisms $H^2(V_i)\to H^2(C_i)$,   $H^2(V_i)\to H^2(C_{i-1})$ are surjective because $V_i$ is ruled over $C_i$ and over $C_{i-1}$. By part (2) of Lemma \ref{uh} one of the rational surfaces, say $V_0$, is not minimal. If $D$ is a smooth rational $-1$-curve
on $V_0$, we have 
$$-1 =(K_{V_0} \cdot D)_{V_0}=   (-C_0 \cdot D)_{V_0}.$$ 
In particular,  the restriction morphism $H^2(V_0)\to H^2(C_0)$ is surjective. Surjectivity of  $\delta_2'$ follows. Thus, the third complex in (\ref{spectral}) has nontrivial cohomology only in the middle degree. As the complex is self-dual, the middle cohomology group must be torsion free.
Consider the fourth complex. Identifying $H^1(C_i)$ with $H^1(E)=:H$, we find that the fourth complex is isomorphic to 
   $$H^{\oplus{m-1}}\rar{\delta_1}H^{\oplus{m}}$$
with the differential given by the formula
$$\delta_1(u_1,\ldots,u_{m-1})=(u_1,u_2-u_1,\ldots, u_{m-1}-u_{m-2},-u_{m-1}).$$
In particular, it has nontrivial cohomology group only in a single degree and this group is isomorphic to $H$.  
The second complex in  (\ref{spectral})  is dual to the fourth one. 
  This completes the proof of degeneration. 
  
  Since the spectral sequence degenerates at $E_2$ and the $E_2$ terms are torsion free it follows that $W_1^{\bZ}=\Coker(\delta_1) \cong H= H^1(E,\bZ)$. Thus $W_1^{\bZ}$ determines the elliptic curve $E$. 
  
  It remains to prove that  $m^2=r_1$.\footnote{This fact is stated without proof in (\cite{fs}).} Indeed, we have the following commutative diagram of abelian groups
 $$
\xymatrixrowsep{.25in}
\xymatrixcolsep{.3in}
\xymatrix{
W_3^\bZ=E_2^{-1,3}\simeq H \ar[rr]^<(.4){\Delta}\ar[ddrrrr]^>(.75){N}|>(.595)\hole
& & H^{\oplus{m}}\ar[rr]^{\delta_3}\ar[dd]_>(.75){Id}& &H^{\oplus{m-1}}\\
& & & & \\
H^{\oplus{m-1}}\ar[rr]^{\delta_1}& &H^{\oplus{m}}\ar[rr]^<(.25){\Sigma}& &H\simeq E_2^{1,0}=W_1^\bZ
}
$$
where $\Delta$ is the diagonal map, $\Sigma$ is the summation map, and $\delta_3$ is given by the formula
$$\delta_3(u_0,\ldots,u_{m-1})=(u_1-u_0,\ldots, u_{m-1}-u_{m-2}).$$
It follows that
$$N=\Sigma \circ \Delta =m \, Id,$$
and thus we have
$$r_1:= |\Coker(W_3^\bZ\rar{N} W_1^\bZ)|= |\Coker(H\rar{m} H)|=m^2.$$
This completes this proof of the Theorem for type II degenerations.
\subsection{Type III Degeneration}\label{t3}
Suppose that $\oX$ is a type III Kulikov degeneration.    In  (\cite{fs},  Prop. 7.1),  Friedman and Scattone proved that the number of triple points of $Y$ is equal to $r_2$.
Then since the Clemens polytope of $Y$ is a triangulation of $S^2$ it follows that the number of double curves in $Y$ is equal to $\displaystyle{\frac{3}{2}r_2}$ and using Euler's formula for triangulizations of a sphere we have that the number of irreducible components  of $Y$
equals $\displaystyle{   \frac{r_2}{2}+2}$.  
We know that each irreducible component $V_i$  of $Y$ is a smooth rational surface and each $C_{j}$ is a smooth rational curve. Thus for each $C_{j}$ we have $[C_{j}]=\bZ+\bZ(-1)$ and since every non-singular rational surface can be obtained by blowing up either the projective plane or a Hirzebruch surface it follows that $[V_i]=\bZ+a_i\bZ(-1)+\bZ(-2)$ for some $a_i\in\bZ_{\ge 0}$. Let $a=\sum_i a_i$.  Then, we have 
\begin{align*}
\int_{X}&=\sum_{ i\in \pi_0(Y^{(0)})}[V_i]-2\sum_{j\in \pi_0(Y^{(1)})}[C_{j}]+3 r_2\bZ\\
&=\left(\frac{r_2}{2}+2\right)\bZ+a\bZ(-1)+\left(\frac{r_2}{2}+2\right)\bZ(-2)-3r_2(\bZ+\bZ(-1))+3r_2\bZ\\
&=\left(\frac{r_2}{2}+2\right)\bZ(-2)+\left(a-3r_2\right)\bZ(-1)+\left(\frac{r_2}{2}+2\right)\bZ
\end{align*} 
Finally, using Proposition \ref{modQ(1)-Q} it follows that 
$$a-3r_2=20-r_2.$$

\begin{rem}\label{r.f.s.} We claim that in notation of \S \ref{t3} the canonical map (\ref{clw0})
\begin{equation}\label{fsb}
W_4^{\bZ}/W_{3}^{\bZ} \rar{\gamma} H_2(Cl(Y)) 
\end{equation} 
is an isomorphism. Indeed, let $x$ be a generator of $W_4^{\bZ}/W_{3}^{\bZ}$,  and let $$\gamma(x)=\sum_{i\in \pi_0(Y^{(2)})} b_i \delta_i, $$ where  $\delta_i$ are $2$-simplices of $Cl(Y)$. 
Then, since $\gamma(x)\in H_2(Cl(Y))$, the boundary of the $2$-dimensional chain  $\sum_{i\in \pi_0(Y^{(2)})} b_i \delta_i$ is $0$.
As the $\delta_i$ form a triangulation of a  compact connected manifold it follows that all the numbers $|b_i|$ are equal one to the other\footnote{Indeed, every $1$-simplex $\epsilon$ of the triangulation has precisely
two $2$-simplices, say $\delta_i$ and $\delta_j$, adjacent to it. Thus, in order to have the coefficient at $\epsilon$ of the boundary of $\gamma(x)$ vanish $|b_i|$ must be equal to $|b_j|$.}. 
If $b$ denotes their common value, we have by the Picard-Lefschetz formula (\ref{picle})
   $$-<Gr N^2 \gamma(x),  \gamma(x)>=  \sum_{i\in \pi_0(Y^{(2)})} b_i^2  = |\pi_0\left(Y^{(2)}\right)| b^2.$$
   The number at the left-hand side of the above formula equals $r_2$. Thus by Friedman-Scattone's result $b=1$ and therefore $\gamma(x)$ is a generator of  $ H_2(Cl(Y)) $.
   
   It follows from a general result of Berkovich explained in the next section that the group $H_2(Cl(Y))$ and morphism (\ref{fsb}) are independent of the choice of a strictly semi-stable model $\oX$. Thus, it is an isomorphism for every such model. 
   
\end{rem}
% \section{Abelian varieties and Kummer K3 surfaces.}\label{a.v.k.s}
 
 %\subsection{The limit mixed Hodge structure.}  With every variation $\cV$ of polarizable pure Hodge structure over a punctured disk $D^*\subset D$  with a fixed coordinate $t$ Schmid associated in \cite{sch} a mixed Hodge structure $\lim \cV= (V_\bZ, W_i^{\bQ}\subset V_\bQ, F^i\subset V_\bC)$. In %addition, the underlying abelian group $V_\bZ$
 %carries an action of the quasi-unipotent monodromy operator $M: V_\bZ \to V_\bZ$.  Moreover, if we write $M_\bQ= U T=TU : V_\bQ \to V_\bQ$,  with unipotent $U$ and semi-simple $T$ both  commuting with $M_\bQ$, then $N=log U$  is a morphism of rational Hodge structures  
 %$$N: \lim \cV \otimes \bQ \to \lim \cV \otimes \bQ(-1)$$
 %and $T$ is an automorphism of   $ \lim \cV \otimes \bQ $ of finite order. While the mixed Hodge structure $\lim \cV$ depends on the choice of coordinate $t$  the pure $\bQ$-Hodge structures $Gr_W^i \lim \cV$ equipped the morphisms $Gr\, T:  Gr_W^i \lim \cV \to Gr_W^i \lim \cV$, $Gr\, T:  Gr_W^i %\lim \cV \to (Gr_W^{i-2} \lim \cV) \otimes \bQ(-1)$ do not.
 
 %If the variation  $\cV$ is of geometric origin then $\lim \cV$ carries a certain additional structure.

 \section{The monodromy pairing.}\label{s.m.p.} 
 Let $K$ be a complete discrete valuation field, and let $\sK$ be the completion of an algebraic closure $\oK$ of $K$.  In (\cite{ber1}), Berkovich  developed a theory of analytic spaces  over  $K$. The underlying topological space $|X^{an}_\sK |$ of the analytification of a scheme $X$ over $K$ has interesting topological invariants (in contrast with the space $X(\sK)$ equipped with the usual 
 topology, which is totally disconnected).  In particular, if $X$ is the generic fiber of a proper strictly semi-stable scheme $\oX$ over $R$ the space $|X^{an}_\sK |$ is homotopy equivalent to the Clemens polytope of the special fiber $Y$ (\cite{ber3}, \S 5). In this section we construct a positive
 pairing on the singular cohomology group  $H^m(|X^{an}_{\sK} |, \bQ)$ that generalizes Grothendieck's monodromy pairing in the case of abelian varieties. Applications to motivic integrals are discussed in the last section.  
  \subsection{Cohomology of the analytic space associated with a smooth scheme.}
  Let $R$ be a  complete discrete valuation domain,  $K$ its fraction field, $k$ the residue field, and let $I\subset G=\mathrm{Gal}(\overline K/K)$ be the inertia subgroup.  
 We denote by $s$ and $\eta$ the closed and generic points of $\mathrm{spec} \,R$ respectively. 
 For a prime number $\ell$ different from $\mathrm{char}\, k$,  we have a canonical surjection (\cite{sga7}, \Rmnum{1},   \S 0.3)
 $$\chi:  I \to \bZ_{\ell}(1)(\overline k).$$ 
 If $\rho: G\to {\mathrm {Aut}}(V)$ is a finite rank $\bZ_{\ell}$-representation of $G$ there is a canonical $G$-homomorphism:
$$N: V\otimes  \bQ_{\ell}(1) \to  V\otimes  \bQ_{\ell},$$
defined as follows. The composition $\Log \circ \rho $ with the $\ell$-adic logarithm ${\mathrm {Aut}}(V)\to {\mathrm {End}}(V\otimes \bQ_{\ell})$ restricted to the inertia subgroup $I$ factors through $\chi$. The map $\bZ_{\ell}(1) \to {\mathrm {End}}(V\otimes \bQ_{\ell})$ yields $N$.

Denote by $\widehat {K}$ the completion of the algebraic closure $\overline K$ with respect to the unique valuation $\overline K^* \to \bQ$  extending the valuation on $K$. For a smooth scheme $X$ of finite type over $K$, 
let $X^{an}_{\sK}$ be the  $\sK$-analytic space associated with $X\otimes _K \sK$  (\cite{ber1}, \S 3.4),  and let  $|X^{an}_\sK |$ be the underlying topological space.  According to (\cite{ber3}, Theorem 9.1;   \cite{hl}, Theorem 13.1.7)  
 $|X^{an}_\sK |$ is a paracompact locally contractible topological space  homotopy equivalent to a finite CW complex. In particular, 
 the singular cohomology groups  $$ \Gamma^m_C(X)= H^m(|X^{an}_\sK |, C) $$  with coefficients in a ring $C$ are finitely generated $C$-modules.  The action of the Galois group $G$ on  $|X^{an}_\sK |$ induces one on $\Gamma^m_C(X)$.
 In (\cite{hl}, Theorem 13.1.8) Hrushovski and Loeser proved that there exists a finite normal extension $K' \supset K$ such that the morphism
 $$H^m(|X^{an}_{K'} |, C) \to H^m(|X^{an}_\sK |, C)= \Gamma^m_C(X)$$
 is an isomorphism\footnote{This result was announced in (\cite{ber3}, Theorem 10.1), however the proof in {\it loc.cit.} is not correct: the assertion on p.82 that a proper  hyper-covering of a scheme $X$ induces a  hyper-covering of the topological space $|X^{an}|$
 is false. Example: take the hyper-covering associated with the $r$-fold \'etale cover $\bG_m\to \bG_m$. If the associated simplicial topological space over $|\bG_m^{an}|$ were a hyper-covering one would get an isomorphism between the cohomology of the contractible space
 $|\bG_m^{an}|$ and the group cohomology $ H^*(\bZ/r\bZ, A)$.  In fact, $ \Gamma^*_A(X)$ is an interesting example of cohomology theory that does not have the \'etale descent property.}.  It follows, that the action of $G$  on $\Gamma^m_C(X)$ 
 factors through a finite quotient $G=\mathrm{Gal}(\oK/K)\epi \mathrm{Gal}(K'/K)$. 
 
 \begin{Th}\label{main2}  For every smooth  variety $X$ and every prime number $\ell\ne \mathrm{char}\, k$, the canonical morphism (\cite{ber2}, Theorem 7.5.4; \cite{ber3}, Theorem 3.2)
  \begin{equation}\label{gammatoet} 
 \gamma: \Gamma^m_{\bZ_{\ell}}(X) \to H^m( X^{an}_\sK, \bZ_{\ell}) \simeq   H^m( X_{\overline K}, \bZ_{\ell})
  \end{equation}
induces an isomorphism of $G$-modules 
\begin{equation}\label{mib}
\Gamma^m_{\bQ_{\ell}}(X) \iso  Im( H^m( X_{\overline K}, \bQ_{\ell})(m)\rar{N^m} H^m( X_{\overline K}, \bQ_{\ell})).
\end{equation}
\end{Th}
 We will write $N^m H^m( X_{\overline K}, \bQ_{\ell})$ for the right-hand side of (\ref{mib}).
 \begin{proof}  Without loss of generality we may assume that $k$ is separably closed and that $X$ is irreducible.
  We first prove the theorem assuming that $X$ is projective and has a strictly semi-stable model  $\oX$ over $R$. In this case, according to a key result of Berkovich (\cite{ber3}, \S 5),  $\Gamma^m_C(X)$  
  is isomorphic to the singular cohomology of the Clemens polytope of the special fiber of $\oX$. On the other hand, we consider the weight filtration $W_i$ on  $H^m( X_{\overline K}, \bQ_{\ell})$ (\cite{rz}, \cite{sa}).
  Interpreting the cohomology of the Clemens polytope  as the weight zero part of $H^m( X_{\overline K}, \bQ_{\ell})$  we find that (\ref{mib}) is equivalent to a special case of Deligne's monodromy conjecture 
  which asserts that, for every integer $0\leq i\leq m$,  the morphism $$N^i: Gr^{m+i}_W H^m( X_{\overline K}, \bQ_{\ell})(i) \to Gr^{m-i}_W H^m( X_{\overline K}, \bQ_{\ell})$$
is an isomorphism.  We prove Deligne's  conjecture for $i=m$ using the method of Steenbrink 
  (who proved it for all $i$ and $k=\bC$).
  To prove the theorem for arbitrary smooth $X$ we show that the functors $\Gamma^m_{\bQ_{\ell}}$ and $N^m H^m$, first, admit transfers for finite morphisms
and, second, take every dominant open embedding $U\mono X$ to an isomorphism.  Finally, we use de Jong's alteration result to complete the proof.

  {\it Step 1.} Assuming that $X$ has a projective strictly semi-stable model   $\oX$ over $R$,  $X\simeq \oX \otimes_R K$. 
 Denote by $D_i$, $i=1, 2, \cdots ,s$ the irreducible components of the special fiber $Y= \oX \otimes k$; 
$$Y^{(q)}= \bigsqcup _{I\subset \{1, \cdots, s\}, |I|=q+1} \bigcap_{i\in I}D_i,$$
and by $\pi_0\left(Y^{(q)}\right)$ the set of connected components of $Y^{(q)}$.
 We have a commutative diagram
 \begin{equation}\label{dcd}
\xymatrix{
  \Gamma^m_{\bQ_{\ell}}(X)   \ar[d]_{\simeq}   \ar[r]_{\gamma }&  H^m( X_{\overline K}, \bQ_{\ell})\\ 
 H^m_{sing} (Cl(Y), \bQ_{\ell}) \simeq E_2^{m,0}(\oX)  \ar[ur]_{\rho} &  \\
}
 \end{equation}
where $E_2^{m,0}(\oX)$  is the weight zero term of the weight spectral sequence converging to $H^m( X_{\overline K}, \bQ_{\ell})$ (\cite{sa}). 
  According to (\cite{na}) the weight spectral sequence degenerates at $E_2$; in particular the morphism $\rho$ is injective.
 Since the range of the weight filtration on $H^m( X_{\overline K}, \bQ_{\ell})$ is at most $2m$ and $N$ shifts the filtration by $2$, we have
  \begin{equation}\label{dconj}
 N^m H^m( X_{\overline K}, \bQ_{\ell}) \subset  Im(\rho).
  \end{equation}
%By Cor. 2.8 of {\it loc. cit.}  $E_2^{m,0}(\oX)$ is isomorphic to the middle cohomology of the complex
 %$$ \cdots \to \bQ_l[\pi_0(Y^{(m-1)})]  \to \bQ_l[\pi_0(Y^{(m)})] \to  \bQ_l[\pi_0(Y^{(m+1)})] \to \cdots $$
 %that also computes $H^*_{sing} (Cl(Y)) \otimes \bQ_l $.
 %According to (\cite{nak}) the weight spectral sequence degenerates in $E_2$ term; in particular $\rho$ induces an isomorphism
 %onto the zero term (the smallest one) of the weight filtration on $H^m( X_{\overline K}, \bQ_l)$.
   Consider the commutative diagram
    \begin{equation}\label{dcd2}
\def\normalbaselines{\baselineskip20pt
\lineskip3pt  \lineskiplimit3pt}
\def\mapright#1{\smash{
\mathop{\to}\limits^{#1}}}
\def\mapdown#1{\Big\downarrow\rlap
{$\vcenter{\hbox{$\scriptstyle#1$}}$}}
\begin{matrix}
 E_2^{-m, 2m}(\oX)(m)&  \twoheadleftarrow  &  H^m( X_{\overline K}, \bQ_{\ell})(m)   \cr
  \mapdown{\overline N^m }& & \mapdown{N^m }  \cr
 E_2^{m, 0}(\oX)  & \mono &  H^m( X_{\overline K}, \bQ_{\ell}),
\end{matrix}
 \end{equation}
The upper horizontal arrow in this diagram is the projection to the weight $2m$ quotient.  We will prove,  following  the method of  (\cite{st1}, \S 5),  that, for every $m$, one has 
  \begin{equation}\label{eq3}
\overline N^m: E_2^{-m, 2m}(\oX)(m) \iso E_2^{m, 0}(\oX). 
\end{equation} 
This trivially holds for $d:= \dim X < m$ because in this case both sides of (\ref{eq3}) equal $0$.  Let us prove   (\ref{eq3}) for $m=d$. Consider the following commutative diagram  (\cite{sa})
 \begin{displaymath}
\xymatrix{
 E^{d-1,0}_1  \ar[r] & H^0\left(Y^{(d)}, \bQ_{\ell}\right)= E^{d,0}_1 \ar[r] & E^{d,0}_2 \to 0 \\
0\to E^{-d,2d}_2(d) \ar[rru]\ar[r] & H^0\left(Y^{(d)}, \bQ_{\ell}\right) = E^{-d,2d}_1(d) \ar[u]|\hole \ar[r] & E^{-d+1, 2d}(d) \\}
\end{displaymath}
where 
$$ E^{d-1,0}_1= H^0\left(Y^{(d-1)}, \bQ_{\ell}\right), \quad E^{-d+1, 2d}(d)=H^2\left(Y^{(d-1)}, \bQ_{\ell}\right)(1),$$ 
the diagonal  morphism is  $\overline N^d$, and the vertical arrow is the identity morphism.  The rows of the above diagram are exact and dual to one another.
In particular, we have a non-degenerate paring
$$<, >:  E^{d,0}_2 \otimes E^{-d,2d}_2(d) \to \bQ_{\ell} $$
that identifies $E^{-d,2d}_2(d)$   with $ H_d (Cl(Y)) \otimes \bQ_{\ell}$.
Next, consider the symmetric form
   \begin{equation}\label{eq2}
   E_2^{-d, 2d}(d)  \otimes   E_2^{-d, 2d}(d)  \to \bQ_{\ell}, \quad x\otimes y  \mapsto <\overline N^d x, y>.
  \end{equation}
 We claim that (\ref{eq2})  is non-degenerate.
 In fact, if  $$x= \sum_{v\in \pi_0\left(Y^{(d)}\right)}  a_v v,  \quad y = \sum_{v\in \pi_0\left(Y^{(d)}\right)} b_v v \in E_2^{-d, 2d}(\oX)(d) \subset \bQ_{\ell}\left[\pi_0\left(Y^{(d)}\right)\right],$$
 we have
 $$  <N^d x, y>= \sum a_v b_v.$$
 Thus (\ref{eq2}) comes by extension of scalars from a positive form
    \begin{equation}\label{dmp}
 H_d (Cl(Y), \bQ)  \otimes H_d (Cl(Y), \bQ) \to \bQ.
    \end{equation}
 This proves that the morphism (\ref{eq3}) is injective; since $\dim\, E_2^{-d, 2d} =\dim \, E_2^{d, 0}$, it must be an isomorphism.

%By Poincar\'e duality  
 %$$E_2^{-d, 2d}(\oX)(d)\simeq ( E_2^{d, 0}(\oX))^* \simeq \Ker(\bQ_l[\pi_0(Y^{(d)})]\to  \bQ_l[\pi_0(Y^{(d-1)})])\simeq  H_d (Cl(Y)) \otimes \bQ_l.$$ 

 Assume that $0<m<d$. Choose an embedding $\oX \mono \bP^N_R$ and a generic  hyperplane section $\oZ= \oX\cap  \bP^{N-d+m}_R $ of dimension $m$; $Z= X\cap  \bP^{N-d+m}_K  $. Then  $\oZ$ is again strictly semi-stable and the embedding $i: \oZ \mono \oX$ induces
 a morphism of spectral sequences  $E_r^{p, q}(\oX) \to E_r^{p, q}(\oZ)$. By the Hard Lefschetz Theorem the composition of the restriction morphism and the Poincar\'e pairing
 $$H^m( X_{\overline K}, \bQ_{\ell})\otimes H^m( X_{\overline K}, \bQ_{\ell})\to  H^m( Z_{\overline K}, \bQ_{\ell})\otimes H^m( Z_{\overline K}, \bQ_{\ell})\to \bQ_{\ell}(-m)$$
 is non-degenerate. The induced isomorphism $H^m( X_{\overline K}, \bQ_{\ell})\to (H^m( X_{\overline K}, \bQ_{\ell}))^*(-m)$ takes $ E_2^{m,0}(\oX)\subset H^m( X_{\overline K}, \bQ_{\ell})$ to $( E_2^{-m, 2m}(\oX))^*\subset  (H^m( X_{\overline K}, \bQ_{\ell}))^*$. Thus
 $$\dim\, E_2^{m,0}(\oX) \leq \dim\, E_2^{-m, 2m}(\oX).$$ 
 Let us show that (\ref{eq3}) is injective. It is enough to check that in  the commutative diagram
       \begin{equation}\label{dcd3}
\def\normalbaselines{\baselineskip20pt
\lineskip3pt  \lineskiplimit3pt}
\def\mapright#1{\smash{
\mathop{\to}\limits^{#1}}}
\def\mapdown#1{\Big\downarrow\rlap
{$\vcenter{\hbox{$\scriptstyle#1$}}$}}
\begin{matrix}
 E_2^{-m, 2m}(\oX)(m)&  \rar{\overline N^m} &   E_2^{m, 0}(\oX)   \cr
  \mapdown{i^*}& & \mapdown{i^* }  \cr
  E_2^{-m, 2m}(\oZ)(m)& \stackrel{\overline N^m}{ \iso} &   E_2^{m, 0}(\oZ).
\end{matrix}
 \end{equation}
the left downward arrow is an injection.  We have
       \begin{equation}\label{dcd4}
\def\normalbaselines{\baselineskip20pt
\lineskip3pt  \lineskiplimit3pt}
\def\mapright#1{\smash{
\mathop{\to}\limits^{#1}}}
\def\mapdown#1{\Big\downarrow\rlap
{$\vcenter{\hbox{$\scriptstyle#1$}}$}}
\begin{matrix}
E_2^{-m, 2m}(\oX)(m)&  \mono & E_1^{-m, 2m}(\oX)(m)= H^0(Y^{(m)}, \bQ_{\ell})      \cr
  \mapdown{i^*}& & \mapdown{i^* }  \cr
  E_2^{-m, 2m}(\oZ)(m)& \rar{} & E_1^{-m, 2m}(\oZ)(m)= H^0(Y^{(m)}\cap \bP^{N-d+m}_k, \bQ_{\ell}).
\end{matrix}
 \end{equation}
In this commutative diagram the upper horizontal arrow is an injection because   the incoming differential $0= E_1^{-m-1, 2m}(\oX)\rar{d_1} E_1^{-m, 2m}(\oX)$ is trivial. The right downward arrow is an injection because $\bP^{N-n+m}_k$ intersects every connected component of  $Y^{(m)}$.
This completes the proof of (\ref{eq3})  and that of (\ref{mib}). 

 {\it Step 2.}  Hrushovski and Loeser proved in (\cite{hl}, Th. 13.1.8) that for every smooth variety $X$ and an open dense subset $U\subset X$ the restriction morphism
   \begin{equation}\label{birgamma}
  \Gamma^m_C(X)\to \Gamma^m_C(U)
   \end{equation}
 is an isomorphism.
 Let us show that the functor at the right-hand side of (\ref{mib}) has the same property: 
  \begin{equation}\label{bir}
  N^m H^m( X_{\overline K}, \bQ_{\ell}) \iso N^m H^m( U_{\overline K}, \bQ_{\ell}).
  \end{equation}
  We first prove (\ref{bir}) in the case when $X$ is the generic fiber of a projective strictly semi-stable pair $(\oX, \oZ=\oZ_f \cup Y)$ over  $R$ ( \cite{dej1}, \S 6.3)  and $j: U\mono X$ is the complement to $Z=\oZ \otimes K$ in $X$.
  Denote by $T$ the special fiber $\oZ_f\otimes k$ of the flat part of  $\oZ$ and by $\oj: Y-T \mono Y$ the embedding. 
  The idea of the following argument (that goes back to  Nakayama (\cite{na})) is the following. 
 When the residue field $k$ is finite (\ref{bir})  can be derived form the Weil conjectures, proven by Deligne,  and the formula  (\ref{eq3}) proven in Step 1 ({\it cf.}  \cite{ber4} p. 672).  In general, the works of Fujiwara, Kato and Nakayama on logarithmic \'etale cohomology (\cite{il2}) imply that $\ell$-adic cohomology groups of $X$ and $U$ depend only
 on the special fibers, $Y$ and $Y-T$, endowed with their natural log structures (that, in turn, are determined by the first infinitesimal neighborhood of $Y$ (resp. $Y-T$) in $\oX$ (resp. $\oX- \oZ_f$)).  Then, a specialization argument enables one to reduce to the finite field case. Let us explain the details.
  
  For a scheme  $S_{log}$ over the log point $(\mathrm{spec}\, k)_{log}$ we denote by $R\tilde \epsilon _*$  the functor from the derived category of $\ell$-adic sheaves on the Kummer  \'etale site, $ S_{log}^{ket} $, to the derived category of $\ell$-adic sheaves on $S$ equipped with an endomorphism of weight $2$ {\it i.e.}, a morphism $N: \cF \to \cF(-1)$  (\cite{il2}, \S 8, p. 308).
   Consider the log structure on the scheme $\oX$ associated with the divisor $Y$, and let
 $Y_{log}= (Y, M_Y)$ be the special fiber with the induced log structure.  According to (\cite{il2}, \S 8, Cor. 8.4.3) the action of the wild inertia $P\subset \mathrm{Gal}(\oK/K)$ on  the complexes of nearby cycles 
 $\Psi \bQ_{\ell}$,  $\Psi Rj_* \bQ_{\ell}$ is trivial.  Therefore we can and we will view the nearby cycles as objects  of the derived category of $\ell$-adic sheaves on  $Y$ endowed with an endomorphism $N$ of weight $2$. Then, we have 
 %The log scheme $(\oX, M_{\oX, Y})$ is smooth over the log disk $(\mathrm{spec}\, R, R-0)$. 
 %Then, it follows from a result of Nakayama that the complexes of nearby cycles $\Psi \bQ_l$,  $\Psi Rj_* \bQ_l$ viewed as complexes of sheaves on $Y$ together with the monodromy
 %action $N: \Psi \bQ_l \to \Psi \bQ_l(-1)$,     $N: \Psi  Rj_*\bQ_l \to \Psi  Rj_* \bQ_l(-1)$,   depend on the log schemes   $\oj: (Y-T)_{log} \mono Y_{log}$ over the log point $(\mathrm{spec}\, k)_{log}$.   Namely,  if $\tilde s \to  (\mathrm{spec}\, k)_{log}$ is a geometric log point over  $\tilde s \to  (\mathrm{spec}\, k)_{log}$,  we have (\cite{il2}, \S 8)
 $$\Psi \bQ_{\ell} \simeq R\tilde \epsilon _*(\bQ_{\ell}),$$
 $$\Psi Rj_* \bQ_{\ell} \simeq R\tilde \epsilon _*(R\oj_* \bQ_{\ell}).$$
% where $R\tilde \epsilon _*$  is a functor from the derived category of l-adic sheaves on the Kummer  \'etale site $ Y_{log}^{ket} \times_s \eta$ to the derived category of l-adic sheaves on $Y$ equipped with an endomorphism of weight $2$ ({\it i.e.}, a map $\cF \to \cF(-1)$).
  We have to prove that the morphism
$$ N^mH^m(Y_{\overline k},   R\tilde \epsilon _*(\bQ_{\ell}) )\to N^mH^m(Y_{ \overline k},  R\tilde \epsilon _*(R\oj_* \bQ_{\ell})   )$$
is an isomorphism. This will follow from a more general fact about log schemes over $(\mathrm{spec}\, k)_{log}$.

Let  $Y_{log}=(Y, M_Y)$ be a fs log scheme over  $(\mathrm{spec}\, k)_{log}$, and let $T\mono Y$ be closed subscheme.  
We will say that $( Y_{log}, T)$ is a standard log strictly semi-stable pair if, for some integers $0\leq a\leq b\leq d$, there is  an isomorphism between $Y_{log}$  and the special fiber the log scheme $\mathrm{spec}\, R[x_0,\cdots x_d] /(x_0\cdots x_a -\pi)$ (with the log structure defined by the divisor $\pi=0$) that takes $T$  to the subscheme given by the equation $x_{a+1}\cdots x_b=0$.    
 We will say that  $( Y_{log}, T)$ is a log strictly semi-stable pair if every point of $Y$ has a Zariski neighborhood $U$ such that $(U_{log}, T\cap U)$ admits a strict \'etale morphism to a standard log strictly semi-stable pair. 
 If this is the case,  every irreducible component $T_i$ of  $T= T_1 \cup \cdots \cup T_n$ with the log structure  induced from $Y$  and  $T_i \cap  (T_1 \cup \cdots  T_{i-1}) \subset T_i$ is again a log strictly semi-stable pair.
 
 Let $( Y_{log}, T)$ be  a proper log strictly semi-stable pair.  In (\cite{na}, \S1), Nakayama constructed the weight spectral sequence $E_r^{pq}$ converging to $H^m(Y\otimes \overline k,   R\tilde \epsilon _*(\bQ_{\ell}) )$ and proved that it degenerates in the $E_2$-terms.
 In particular, for every integer $m$, the canonical morphism 
 $$ H^m_{sing} (Cl(Y)) \otimes \bQ_{\ell} \simeq E_2^{m,0} \to  H^m(Y\otimes \overline k,   R\tilde \epsilon _*(\bQ_{\ell}) )$$
 is an embedding.
  \begin{lm}\label{brr}
 For every proper log strictly semi-stable pair $( Y_{log}, T)$ the composition
  \begin{equation}\label{logbir}
H^m_{sing} (Cl(Y_{\overline k}), \bQ_{\ell})\mono H^m(Y_{\overline k},   R\tilde \epsilon _*(\bQ_{\ell}) ) \to H^m(Y_{\overline k},  R\tilde \epsilon _*(R\oj_* \bQ_{\ell})   )
   \end{equation}
is a monomorphism whose image contains $ N^mH^m(Y_{\overline k},  R\tilde \epsilon _*(R\oj_* \bQ_{\ell})   ) $.
  \end{lm}
 \begin{proof} The specialization argument of Nakayama (\cite{na}) reduces the statement to the case when $k$ is a finite field; in the rest of the proof we will be assuming that this is the case.
 The vector spaces appearing in (\ref{logbir}) carry an action of the Galois group $\mathrm{Gal}( \overline k/k)$. 
 Let us look at the action  of the Frobenius element  $Fr\in \mathrm{Gal}( \overline k/k)$.   For a finite-dimensional $\ell$-adic representation $V$ of $\mathrm{Gal}( \overline k/k)$ we denote by $V_0$ 
 the largest invariant subspace of $V$ such that all the eigenvalues of $Fr$ on $V_0$ are roots of unity.  Looking at the weight spectral sequence we see that 
 $$  E_2^{m,0}   =(H^m(Y_{\overline k},   R\tilde \epsilon _*(\bQ_{\ell}) )_0.$$
 Thus, to prove the lemma it suffices to show the following:
 \begin{enumerate}[(a)]
\item
  \begin{equation}\label{berw}
\left(H^m(Y_{\overline k},   R\tilde \epsilon _*(\bQ_{\ell}) )\right)_0\iso \left(H^m(Y_{\overline k},  R\tilde \epsilon _*(R\oj_* \bQ_{\ell})   )\right)_0.
 \end{equation}
 \item
 The eigenvalues of $Fr$ acting on $H^m(Y_{\overline k},  R\tilde \epsilon _*(R\oj_* \bQ_{\ell})   )$ are Weil numbers of weights from $0$ to $2m$.
  \end{enumerate}
   Arguing by induction on $d=\dim Y$ we assume that the above assertions hold for log strictly semi-stable pairs of dimension less then $d$. Let $T_1,\cdots T_n$ be irreducible components of $T$,
  let $Y_j$ be the complement to $ \bigcup_{i\leq j} T_i$ in $Y$.  Consider the Gysin exact sequence
  $$\cdots \to H^{m-2}\left((T_{j+1}\cap Y_j) \otimes \overline k  ,   R\tilde \epsilon _*(\bQ_{\ell}) \right)(-1)\to H^{m}(Y_j  \otimes \overline k,   R\tilde \epsilon _*(\bQ_{\ell}) )$$
  $$ \to  H^{m}(Y_{j+1}  \otimes \overline k,   R\tilde \epsilon _*(\bQ_{\ell}) )\to H^{m-1}\left((T_{j+1}\cap Y_j)  \otimes \overline k,   R\tilde \epsilon _*(\bQ_{\ell}) \right)(-1) \to \cdots $$
  By our induction assumption the boundary terms of the sequence have weights between $2$ and $2m$.  Induction on $j$  proves the first claim (\ref{berw}). 
  The second claim also follows from the above and from the fact that  $H^m(Y\otimes \overline k,   R\tilde \epsilon _*(\bQ_{\ell}) )$ has weights between  $0$ and $2m$.
   \end{proof}
   
 As we know from Step 1, for a projective strictly semi-stable  scheme $\oX$ over $R$, we have
$$N^m H^m( X_{\overline K}, \bQ_{\ell})  \iso E_2^{m,0}.$$
This together with the Lemma \ref{brr}  complete the proof of (\ref{bir}) for strictly semi-stable pairs.

Before going further,  recall that, for every  generically finite surjective morphism $f: X' \to X$ of smooth connected varieties, the induced map
$$f^*:  H^m( X_{\overline K}, \bQ_{\ell}) \to H^m( X'_{\overline K}, \bQ_{\ell})$$
is injective. In fact,  the canonical isomorphism $ \bQ_{\ell} \iso Rf^! \bQ_{\ell}$ defines by adjunction a morphism
$$Rf_* \bQ_{\ell} \iso Rf_!  \bQ_{\ell} \to \bQ_{\ell}.$$ 
In turn, the latter yields the transfer morphism
$$f_*:  H^m( X'_{\overline K}, \bQ_{\ell}) \to H^m( X_{\overline K}, \bQ_{\ell})$$
such that the composition $f_*f^*$ equals  multiplication by the degree of $f$ over the generic point.

Let us return to the proof of (\ref{bir}). Without loss of generality we may assume that $X$ is connected.
Then, by de Jong's result ( \cite{dej1}, \S 6.3)   
we can find a proper generically finite surjective morphism $f: X' \to X$ such that $X'$ is an open subscheme of a connected projective strictly semi-stable scheme  $\oX'$ over  a finite extension $R'\supset R$ and such that  $(\oX', \oX' -X')$ is a strictly semi-stable pair.
Applying de Jong's result once again, we find a proper generically finite surjective morphism $g: \oX^{\prime \prime} \to \oX'$, with connected $\oX^{\prime \prime}$, such that $(\oX^{\prime \prime} , \oX^{\prime \prime} - (fg)^{-1}(U))$ is a projective strictly semi-stable pair over some $R^{\prime \prime} \supset R'$. Diagram:
       \begin{equation}\label{dcd49}
\def\normalbaselines{\baselineskip20pt
\lineskip3pt  \lineskiplimit3pt}
\def\mapright#1{\smash{
\mathop{\to}\limits^{#1}}}
\def\mapdown#1{\Big\downarrow\rlap
{$\vcenter{\hbox{$\scriptstyle#1$}}$}}
\begin{matrix}
(fg)^{-1}(U) &  \mono & X^{\prime \prime} &\mono & \oX^{\prime \prime}    \cr
  \mapdown{}& & \mapdown{ } &&  \mapdown{g }    \cr
f^{-1}(U) &  \mono & X^{\prime} &\mono & \oX^{\prime}      \cr
  \mapdown{}& & \mapdown{f } &&    \cr
    U&\mono & X &&   
\end{matrix}
 \end{equation}
We know that (\ref{bir}) is true for the embeddings $X^{\prime} \mono  \oX^{\prime} \otimes K$ and $g^{-1}f^{-1}(U)   \mono  \oX^{\prime \prime} \otimes K$.\footnote{Indeed, $(\oX', \oX' -X')$ is a strictly semi-stable pair over $R'$. Therefore, we have
 $N^m H^m( \oX' \times _{R'} \oK, \bQ_{\ell}) \iso N^m H^m( X' \times _{R'} \oK, \bQ_{\ell})$. This implies that the morphism $N^m H^m( \oX' \times _{R} \oK, \bQ_{\ell}) \to N^m H^m( X' \times _{R} \oK, \bQ_{\ell})$ is an isomorphism as well.}
Define a morphism  $u: N^m H^m( U_{\overline K}, \bQ_{\ell})\to N^m H^m( X_{\overline K}, \bQ_{\ell}) $
to be the composition
$$N^m H^m( U_{\overline K}, \bQ_{\ell}) \rar{(fg)^*} N^m H^m\left( (fg)^{-1}(U)_{\overline K}, \bQ_{\ell}\right) \simeq N^m H^m(  \oX^{\prime \prime} \otimes \oK, \bQ_{\ell})$$
$$\rar{g_*}N^m H^m(\oX^{\prime} \otimes \oK, \bQ_{\ell})\rar{Res}N^m H^m(X^{\prime} \otimes \oK, \bQ_{\ell})
\rar{f_*}N^m H^m(X \otimes \oK, \bQ_{\ell}).$$
An easy diagram chase shows that $u$ divided by the degree of the morphism  $fg$ over the generic point is the two-sided inverse to the restriction morphism (\ref{bir}).
 
\textit{Step 3.} Let $f: U' \to U$ be a finite surjective morphism of connected smooth varieties. Assume that the corresponding extension $\mathrm{Rat}(X) \subset \mathrm{Rat}(X')$  of the field of rational functions is normal and let $G$ be its Galois group. Then, the pullback morphism $f^*$ induces an isomrphism  
 $$N^m H^m(U_\oK, \bQ_{\ell}) \iso  (N^m H^m(U^{\prime}_ \oK, \bQ_{\ell}))^G.$$
Let us show the functor $\Gamma^m_\bQ$ has the same property:
   \begin{equation}\label{descent}
    \Gamma^m_\bQ (U) \iso  (\Gamma^m_\bQ (U'))^G. 
       \end{equation}
Indeed, by (\cite{ber2}, Prop. 4.2.4),  the cohomology of the topological space $|U^{an}_\sK |$ with rational coefficients coincides
 with the \'etale cohomology of the analytic space $U^{an}_\sK$ with coefficients in $\bQ$. Next, since the functor of $G$-invariants is exact in any $\bQ$-linear abelian category, we have
  $$(H^m_{et}(U^{\prime an}_\sK, \bQ))^G  \simeq H^m_{et}(U^{an}_\sK, (f_* \bQ)^G). $$
  We complete the proof of (\ref{descent}) by showing that the canonical morphism $\bQ \to (f_* \bQ)^G$ is an isomorphism.
  In fact, the weak base change theorem (\cite{ber2} Th. 5.3.1) reduces the statement to the case when $U_\sK$ is a single point. In this case $G$ acts transitively on points of $U'_\sK$ and our assertion follows. 
  
\textit{Step 4.} Now we can complete the proof of Theorem \ref{main2}.  We may assume that $X$ is connected. Then, by (\cite{dej2}, Th 5.9), there exists a proper generically finite surjective morphism $f: X' \to X$ such that  the field extension  $\mathrm{Rat}(X) \subset \mathrm{Rat}(X')$ is normal,  $X'$ is an open subscheme of a connected projective strictly semi-stable scheme  $\oX'$ over  a finite extension $R'\supset R$.  Let $U$ be an open dense subset of $X$ over which $f$ is finite.
By the result of Step 1 the Theorem is true for $\oX'$.\footnote{Indeed, the result of Step 1 implies that the morphism $H^m( |X^{\prime an} \times_{K'} \sK|, \bQ_{\ell}) \to   N^m H^m( X'\times_{K'} \oK , \bQ_{\ell})$ is an isomorphism. This implies
that  $H^m( |X^{\prime an} \times_{K} \sK|, \bQ_{\ell}) \to   N^m H^m( X'\times_{K} \oK , \bQ_{\ell})$  is also an isomorphism.}  Then, by  Step 2 it is true for $f^{-1}(U)$ and thus, by Step 3, for $U$. Applying the result of Step 2 once again we complete the proof of Theorem \ref{main2}.

 \end{proof}

% \begin{rem} According to (\cite{ber3}) is a birational invariant of $X/K$.
 %\end{rem}
 \begin{rem}\label{w.0.a} 
 The groups  $ \Gamma^*_\bZ(X)$ are related to the weight zero part of motivic  vanishing cycles $\Psi(X)\in DM_{gm}^{eff}(k)$ of $X$ (\cite{a1}, \cite{a2}).  Namely, if $\mathrm{char}\, k=0$, 
 %\footnote{We do not know how to define $\Psi(X)$ for an arbitrary d.v.r. $R$. }
  one has
    $$ \Gamma^m_\bZ(X) \simeq Hom_{DM_{gm}^{eff}(\overline k)}(\Psi(X), \bZ[m]).  $$  
 \end{rem} 
\begin{rem}\label{w.0.l.h}  Assume that $K=\bC((t))$.  For every smooth projective $X/K$ there is a canonical morphism ({\it cf.}  \cite{ber5}, Theorem 5.1)
\begin{equation}\label{berhodge} 
 \Gamma^m_\bZ (X) \to  W_0^{\bQ} \cap  H^m(\lim X, \bZ)  
 \end{equation}
 that induces an isomorphism modulo torsion 
   \begin{equation}\label{hodgeber} 
\Gamma^m_\bQ (X)  \simeq W_0^{\bQ} H^m(\lim X).
   \end{equation}
  Morphism  (\ref{berhodge}) can be constructed as follows. Pick a finite extension $K'\supset K$ and strictly semi-stable model $\oX_{R'}$ of  $X_{K'}= X\otimes _K K' $ over the integral closure $R'$ of $R$ in $K'$. Then
   (\ref{berhodge}) is defined to be the composition
   \begin{eqnarray}\label{berhodgecon} 
   \Gamma^m_\bZ (X) \simeq \Gamma^m_\bZ (X_{K'})\iso  H^m(Cl(Y)) \\
  \nonumber \to  W_0^{\bQ} \cap  H^m(\lim X_{K'}, \bZ)\simeq W_0^{\bQ} \cap  H^m(\lim X, \bZ), 
    \end{eqnarray} 
  where $Y$ is the special fiber of $\oX_{R'}$ and the map $H^m(Cl(Y)) \to  W_0^{\bQ} \cap  H^m(\lim X_{K'}, \bZ)$  comes from the weight spectral sequence (see \S \ref{l.h.s}). As the weight spectral sequence with rational coefficients degenerates at $E_2$ terms 
  the above composition is an isomorphism up to torsion. The composition of  (\ref{berhodgecon}) with the embedding  $W_0^{\bQ} \cap  H^m(\lim X, \bZ) \mono H^m( X_{\overline K}, \bZ_{\ell})$ equals
  the canonical morphism $\Gamma^m_\bZ (X) \to H^m( X_{\overline K}, \bZ_{\ell})$ from Theorem \ref{main2}. Thus, the morphism $ \Gamma^m_\bZ (X) \to  W_0^{\bQ} \cap  H^m(\lim X, \bZ)  $ induced
  by (\ref{berhodgecon}) is independent of the choice of $K'$ and $\oX_{R'}$. 
     
 In general, morphism (\ref{berhodge}) is not bijective. 
 \end{rem}

We conjecture that for every smooth proper variety  $X$ over $K$, one has 
\begin{equation}\label{mlemma}  
 \dim_\bQ  \Gamma^m_\bQ (X)  \leq  \dim _K H^m(X, \cO_X).
  \end{equation}
Conjecture  (\ref{mlemma}) is motivated by the following result. 
\begin{pr} The inequality  (\ref{mlemma}) is true if either of the following conditions holds.
\begin{enumerate}[(a)]
 \item  $\mathrm{char}\, k=0.$
\item  $K$ is a finite extension of $\bQ_p$.
 \end{enumerate} 
 \end{pr}
 \begin{proof} When proving the first part of the Proposition, we may assume that $R=\bC[[t]]$ and $X$ is the generic fiber of a strictly semi-stable scheme $\oX$ over $R$ (\cite{hl}, Theorem 13.1.8). In this case, we have
  $$\Gamma^m_\bQ(X) \simeq H^m_{Zar}(Y, \bQ)\simeq W_0^{\bQ} H^m(\lim X).$$
  where $Y$ is the special fiber of $\oX$. The first part of the Proposition now follows from the inequality $ \dim_\bQ   W_0^{\bQ} \leq \dim_\bC   F^0/F^1= \dim _K H^m(X, \cO_X)$. For the second part, recall from (\cite{ber2}, Theorem 1.1)
  that $\Gamma^m_K(X) $ is isomorphic to the subspace of the p-adic cohomology $H^m(X_{\oK}, \bQ_p) \otimes _{\bQ_p} K$ that consists of smooth vectors {\it i.e.}, vectors whose stabilizer in $G$ is open. Thus,
  $$\dim _{K} \Gamma^m_K(X) \leq \dim_{K} (H^m(X_{\oK}, \bQ_p)\otimes _{\bQ_p} \bC_p)^G= \dim _K H^m(X, \cO_X).$$
  The last equality follows from the Hodge-Tate decomposition proven by Faltings (\cite{fa}). 
 \end{proof}

\subsection{The monodromy pairing.}\label{m.p.def}
Let $X$ be a  smooth variety over a complete discrete valuation  field $K$ and  $d= \dim X$.  
In this subsection we define a canonical positive symmetric form  (that we shall call the monodromy pairing)
 \begin{equation}\label{mp} 
 (\cdot, \cdot ):  \Gamma^d_\bQ (X)   \otimes   \Gamma^d_\bQ (X)   \to \bQ.
 \end{equation}
 %where $\Gamma^d_{tor}(X)$   denotes the torsion part of $\Gamma^d(X)$.  
 The group $\Gamma^d_\bZ(X)$ as well as the monodromy pairing depends only on the class of $X$ modulo birational equivalence.

 First, we define a pairing   
 $$ (\cdot, \cdot )_{\ell}: N^d H^d( X_{\overline K}, \bQ_{\ell})\otimes  N^d H^d( X_{\overline K}, \bQ_{\ell}) \to \bQ_{\ell}.$$
  By (\cite{dej1}, Th 4.1, Rem. 4.2), there exists a proper generically finite surjective morphism $f: X' \to X$ such that  $X'$ is an open subscheme of a smooth projective  variety  $\tilde X' $ over  a finite extension $K'\supset K$.
 Let $r$ be the degree of $f$ over the generic point. 
 Consider the morphism
 $$\of^*: N^d H^d( X_{\overline K}, \bQ_{\ell}) \rar{f^*}  N^d H^d( X' _\oK, \bQ_{\ell}) \stackrel{ \sim} {\leftarrow} N^d H^d( \tilde X'_\oK, \bQ_{\ell}).
 \footnote{We write $X'_\oK$ for the fiber product of $X'$ and $\mathrm{spec} \, \oK$  over $\mathrm{spec} \, K$.}$$
 The left arrow is an isomorphism by (\ref{bir}). Given $x, y\in N^d H^d( X_{\overline K}, \bQ_{\ell})$ we set
  \begin{equation}\label{eqmpladic}
  (x, y )_{\ell}= \frac{(-1)^{\frac{d(d-1)}{2}}}{r} <\of^*(x), \of^*(y')>,
  \end{equation}
 where $y' \in  H^d( X_{\overline K}, \bQ_{\ell})(d) $ is an element such that $N^dy'=y$ and $$<,>: H^d( \tilde X'_\oK, \bQ_{\ell})\otimes H^d( \tilde X' _\oK, \bQ_{\ell})(d) \to \bQ_{\ell}$$
  is the Poincar\'e pairing.
 Let us check that $(\cdot, \cdot )_{\ell}$ is well defined. Indeed, if $y^{\prime \prime}$ is another element such that $N^dy^{\prime \prime}=y$, then
 $$<\of^*(x), \of^*(y' -y^{\prime \prime})>= <N^d\of^*(x'), \of^*(y' -y^{\prime \prime})>$$
 $$=(-1)^d  <\of^*(x'), N^d\of^*(y' -y^{\prime \prime})>=0.$$
 The independence of the choice of $X'$, $ f $ and $\tilde X'$ follows from the fact that given another such triple  $X^{\prime \prime}$, $ g$ and $\tilde X^{\prime \prime} $ we can find 
 a smooth projective scheme over some finite extension of $K$ that admits proper generically finite surjective morphisms to both   $\tilde X'$ and $\tilde X^{\prime \prime} $.
 
   Let us also remark that  the pairing $(\cdot, \cdot )_{\ell}$ is symmetric.
   \begin{Th}\label{main3}  
 For every smooth connected variety $X$ of dimension $d$,
the restriction (\ref{mp}) of the pairing  $ (\cdot, \cdot )_{\ell}$ to the subspace 
   $$\Gamma^d_\bQ(X) \mono N^m H^m( X_{\overline K}, \bQ_{\ell})$$
    takes  values in $\bQ$ and is independent of ${\ell}\ne \mathrm{char} \, k$.
  Moreover, the pairing (\ref{mp}) is positively defined (and, in particular,  non-degenerate).
 \end{Th}
 \begin{proof}  Thanks to the birational invariance property of  $\Gamma^d(X)$ (\ref{birgamma})  and de Jong's semi-stable reduction theorem ( \cite{dej1}, \S 6.3) it is enough to prove the theorem in the case when $X$ is the generic fiber of a strictly semi-stable projective scheme $\oX$ over $R$.
In this case, using (\ref{dcd}) we have a canonical isomorphism $\Gamma^d_\bQ (X)\simeq  H^d (Cl(Y), \bQ)$ that identifies, by the Picard-Lefschetz formula ({\it cf.} (\ref{picle})), the pairing $(\cdot, \cdot )_{\ell}$ restricted to $\Gamma^d_\bQ(X)$ with the dual of the pairing (\ref{dmp}).
\end{proof}

 \begin{rem}\label{m.p.p.v.}  The construction of the monodromy pairing can be generalized as follows.  For a pair $(X, \mu)$, where $X$ is a smooth projective variety over $K$ and  $\mu \in H^2(X, \bQ_{\ell}(1)) $ is the class of an ample line bundle over $X$,  and an integer $ m\leq d$,  we define  
 a positive symmetric form
  \begin{equation}\label{mppol} 
 (\cdot, \cdot )_{\mu}: \Gamma^m_\bQ(X)   \otimes \Gamma^m_\bQ(X)  \to \bQ
 \end{equation}
  to be the composition  
$$\Gamma^m_\bQ(X)   \otimes \Gamma^m_\bQ(X) \to N^m H^m( X_{\overline K}, \bQ_{\ell})\otimes N^m H^m( X_{\overline K}, \bQ_{\ell}) \rar{  (\cdot, \cdot )_{\ell, \mu}} \bQ_{\ell}, $$
 where $(x, N^m y' )_{\ell, \mu}= (-1)^{\frac{d(d-1)}{2}}  <x, y' \mu^{d-m}>$.  Let us prove  that (\ref{mppol}) is independent of $\ell$ and positive. Without loss of generality, we may assume that $\mu$ is the class of very ample line bundle $L$. Let $X\mono \bP^N_K$ be the corresponding
 embedding, and let  $Z= X\cap  \bP^{N-d+m}_K\stackrel{i}{\mono} X$ be a generic  hyperplane section of dimension $m$. Then,    $(\cdot, \cdot )_{\mu}$ equals the composition
 $$\Gamma^m_\bQ(X)   \otimes \Gamma^m_\bQ(X) \rar{i^*\otimes i^*} \Gamma^m_\bQ(Z)   \otimes \Gamma^m_\bQ(Z) \rar{ (\cdot, \cdot )} \bQ.$$
 By Theorem \ref{main2} and the Hard Lefschetz Theorem the restriction morphism $i^*: \Gamma^m_\bQ(X)\to \Gamma^m_\bQ(Z)$ is injective. Our claim follows from Theorem \ref{main3}.
 \end{rem}
 \begin{rem}\label{m.p.via.l.h}  Assume that $K=\bC((t))$.  For a smooth projective $d$-dimensional scheme $X$ over $K$  the isomorphism 
$$\Gamma^d_\bQ (X)  \simeq W_0^{\bQ} H^d(\lim X)$$
from Remark \ref{w.0.l.h} carries the monodromy pairing on $\Gamma^d_\bQ (X)$
to the pairing 
$$(\cdot, \cdot):  W_0^{\bQ} H^d(\lim X) \otimes W_0^{\bQ} H^d(\lim X) \to \bQ$$
defined by the formula ({\it cf.} (\ref{eqmpladic}))      
$$(x, y )= (-1)^{\frac{d(d-1)}{2}} <x, y'>,$$
where $x\in W_0^{\bQ}$, $y'\in W_{2d}^{\bQ}/W_{2d-1}^{\bQ}$ is such that $Gr N^d (y')= y$, and  $<\cdot ,  \cdot>: W_0^{\bQ} \otimes W_{2d}^{\bQ}/W_{2d-1}^{\bQ} \to \bQ$ denotes the Poincar\'e pairing.
\end{rem}

\begin{example} Let $A$ be a $d$-dimensional abelian variety over $K$  with semi-stable reduction.
According to (\cite{ber1}, \S 6.5), after replacing $K$ by a finite unramified extension, we can represent the analytic space $A^{an}$ as the quotient of $G^{an}$ by  $\Lambda$, where $G^{an}$ is the analytic group associated with a semi-abelian variety $0\to \overline T\to \overline G \to \overline B\to 0$ over $R$ and $\Lambda  \stackrel{u}{\mono} \overline G(K)$ a lattice.
 Moreover,  the map $|G^{an}_\sK| \to |A^{an}_\sK|$ exhibits $|G^{an}_\sK|$ as a universal cover of $|A^{an}_\sK|$. In particular, $\Gamma_m(A):=H_m(|A^{an}_\sK|)  \simeq \bigwedge^m \Lambda$.   A polarization, $\mu$, of $A$ defines an isogeny
 $\mu _*:  \Lambda \to \Xi $,
  where $\Xi$ is the group of characters of $\overline G$. Using (\cite{c}, Theorem 2.1), we see that the pairing 
  $$\Gamma_1(A)   \otimes \Gamma_1(A)  \to \bQ$$
derived from (\ref{mppol}) equals the pullback of Grothendieck's monodromy pairing
  $$\Lambda \otimes \Xi \rar{u\otimes Id} \overline G(K)/\overline G(R) \otimes \Xi \to \Xi^* \otimes \Xi \to \bZ$$
 via $Id \otimes \mu_*: \Lambda \otimes \Lambda  \to \Lambda \otimes \Xi $, divided by the degree $\mu^d\in \bZ$ of the polarization. 
 \end{example}

\subsection{A birational invariant.}\label{bir.inv}
Let $X$ be a smooth connected variety over a complete discrete valuation  field $K$ and $d= \dim X$.  Assume that $\Gamma^d(X)_\bQ \ne 0$. Let $\mathrm{Disc}(\cdot, \cdot )\in \bQ^*$ be the discriminant of  the monodromy pairing (\ref{mp}) relative to the lattice
$\Gamma^d_\bZ(X)/\Gamma^d_\bZ(X)_{tor} \subset \Gamma^d_\bQ(X)$, 
 and let
\begin{equation}\label{eq.bir.inv}
r_d(X, K)= \frac{1}{\mathrm{Disc}(\cdot, \cdot )}.
\end{equation}
Since the group $\Gamma^d_\bZ(X) $ and the monodromy pairing (\ref{mp}) are birational invariants of $X$ so is the number  $r_d(X, K)$.
If $K\subset K'$ is a finite extension of ramification index $e$, we have 
$$r_d(X\otimes K', K')= e^{d\, \dim \, \Gamma^d_\bQ(X)} r_d(X, K).$$

In the remaining part of this section we shall relate the invariant $r_d(X, K)$ defined here to the one introduced in \S  \ref{MR} for K3 surfaces over $\bC((t))$.
\begin{pr}\label{refrequest} Let $X$ be a smooth projective K3 surface over $K=\bC((t))$ and let  
  $H^2(\lim X)$ be the corresponding limit mixed Hodge structure (see  \S \ref{l.h.s}).   Set $ W_i^{\bZ}: = W_i^{\bQ}\cap H^2(\lim X,\bZ)$. 
  Assume that the monodromy acts on $H^2(\lim X,\bZ)$ by a unipotent operator and let
  $N:  H^2(\lim X, \bZ)  \to H^2(\lim X, \bZ)$ be its logarithm (which is integral by  (\cite{fs}, Prop. 1.2)). 
  Then 
  \begin{enumerate}[(a)]
\item The topological space $|X^{an}_\sK| $  is contractible unless $N^2\ne 0$. If $N^2\ne 0$ the space $|X^{an}_\sK| $ is homotopy equivalent to a $2$-dimensional sphere and the canonical
map (see Remark \ref{w.0.l.h})
 \begin{equation}\label{satk3}
\Gamma^2_\bZ(X) \to W_0^\bZ 
\end{equation}
 is an isomorphism.
 \item  Assume that $N^2\ne 0$. Then the number $r_2(X, \bC((t)))$ defined by (\ref{eq.bir.inv}) is equal to the order of the following group 
     \begin{equation}\label{m.p.repeat}
\Coker\, (W_{4}^{\bZ}/W_{3}^{\bZ}  \rar{Gr \, N^2}  W_{0}^{\bZ}). 
\end{equation}
\end{enumerate}
  \end{pr}
\begin{proof}   It is enough to prove the proposition in the case where $X$ is the restriction of a strictly semi-stable family over a smooth curve. Indeed, at the
expense of a finite extension of $K$ we may choose a strictly semi-stable model $\oX$ for $X$.  The  space $|X^{an}_\sK| $ is homotopy equivalent to the Clemens polytope of the special fiber $Y$ of $\oX$ (\cite{ber3}, \S 5). Applying  Proposition \ref{artinappr} to $\oX$  we find a proper strictly semi-stable family  $\overline \cX$ over a smooth  pointed curve $a\in \oC$, whose fiber over the first infinitesimal neighborhood of point $a$ is isomorphic
to  $\oX\otimes _R R/t^2$ and whose generic fiber is a K3 surface. 
As the limit mixed Hodge structure of a strictly semi-stable scheme depends only
 on the first infinitesimal neighborhood of special fiber the validity of the proposition  does not change  if we replace 
$X$  by $\cX \times_C \mathrm{spec}\, K'$. 
 
   Thus, we may assume that $X$ has a Kulikov model over $R=\bC[[t]]$ (see \S \ref{ss.kulikovmodel}).  If $\oX$ is a Kulikov model, then the Clemens polytope $Cl(Y)$ of the special fiber of $\oX$ is homeomorphic 
   to a point or to an interval for type {\rm I} or {\rm II} degenerations and it is homeomorphic to a $2$-dimensional sphere for type {\rm III} degenerations. This proves the first part of the proposition except for the claim
   that morphism (\ref{satk3}) is an isomorphism. Using Berkovich's result (\cite{ber3}, \S 5), the latter is equivalent the following assertion: the canonical morphism
   $$H^2(Cl(Y)) \to W_0^\bZ$$
   coming from the weight spectral sequence (see Theorem \ref{lhs}) is an isomorphism. In fact, the (equivalent) dual statement,
 $$W_4^{\bZ}/W_{3}^{\bZ} \iso H_2(Cl(Y))$$ 
is proven (using a deep result  of Friedman-Scattone (\cite{fs})) in Remark  \ref{r.f.s.}. This completes the proof of the first part of the proposition.

Part (b) of the proposition follows from the fact that  (\ref{satk3}) is an isomorphism and Remark \ref{m.p.via.l.h}.
\end{proof}

 \section{Motivic integral of maximally degenerate K3 surfaces over non-archimedean fields.}\label{m.i.m.g.k3}
% In this section we state a conjectural formula for the motivic integral of maximally degenerated K3 surface over an arbitrary  non-archimedean field $K$.  First, using
 %the  non-archimedean analytic geometry in the sense of Berkovich (\cite{ber1}) we define a certain cohomological invariant of an arbitrary variety $X$ over $K$.  This invariant should viewed as the weight
 %zero part of the not yet constructed limit motive of $X$. 

 \subsection{The formula.} 
  Throughout this section $R$ denotes a  complete discrete valuation ring with fraction field $K$ and perfect residue field $k$. 
 We shall say that a smooth projective $d$-dimensional Calabi-Yau variety $X$ over $K$ is maximally degenerate if  $\Gamma^d_\bQ(X) \ne 0$.
  According to Theorem \ref{main2}, $X$ is   maximally degenerate if and only if for some (and, therefore, for every) prime $\ell\ne \mathrm{char} \, k$ the map
  $$H^d( X_{\overline K}, \bQ_{\ell})(m)\rar{N^d} H^d( X_{\overline K}, \bQ_{\ell})$$
  is not zero.
\begin{conj}\label{mainc} Let $X$ be a smooth projective maximally degenerate K3 surface over $K$.  Then
\begin{enumerate}[(a)]
\item   The topological space $|X^{an}_\sK| $  is homotopy equivalent to a $2$-dimensional sphere. In particular,   the group $\Gamma^m_\bZ(X)$ is trivial for  $m\ne 0,2$ and isomorphic to $\bZ$ otherwise.
\item  For every $\ell \ne \mathrm{char}\, k$  the lattice  
$$\bZ_{\ell} \simeq \Gamma^2_{\bZ_{\ell}}(X) \stackrel{(\ref{gammatoet})}{\mono} H^2( X_{\overline K}, \bZ_{\ell}) $$
is saturated {\it i.e.}, 
$$  \Gamma^2_{\bZ_{\ell}}(X) = (\Gamma^2_{\bQ_{\ell}}(X))     \cap   H^2( X_{\overline K}, \bZ_{\ell}).$$ 
\item There exists a finite extension $K'\supset K$ such that for every finite extension $L\supset K'$ of ramification index $e$
\begin{eqnarray}\label{conjfor}
\int_{X_L}= \left(\frac{e^2 r_{2}(X,K)}{2}+2\right)\bQ(0)+(20-e^2 r_{2}(X,K))\bQ(-1)\qquad \quad&\\
\nonumber +\left(\frac{e^2 r_{2}(X,K)}{2}+2\right)\bQ(-2).&
\end{eqnarray}
\end{enumerate}
 
\end{conj}

\begin{rem} According to the first part of the conjecture,  for every prime $\ell \ne \mathrm{char}\, k$,  the $\ell$-primary  factor of $r_2(X, K)$ has the following cohomological interpretation. 
If $r_2(X, K)=  \ell^{a_{\ell}} r'$ and $(r', \ell)=1$, then
$$a_{\ell}= -v_{\ell} ({\mathrm {Disc}}(\cdot, \cdot )_{\ell})$$ 
where $ {\mathrm {Disc}}(\cdot, \cdot )_{\ell}$ is the discriminant of  the $\ell$-adic monodromy pairing (\ref{dcd}) restricted to 
$  \mathrm{Im}( N^2) \cap  H^2_{et}( X_{\overline K}, \bZ_{\ell})$ 
 and $v_{\ell}: \bQ_{\ell}^*/\bZ_{\ell}^* \to  \bZ$ is the valuation morphism.  
\end{rem}

\begin{rem} According to Theorem \ref{main} and Proposition \ref{refrequest},  Conjecture  \ref{mainc} is true for $K=\bC((t))$. Thus, it true for every $K$ of equicharacteristic $0$ (\cite{hl}, Theorem 13.1.8).
\end{rem}

 \subsection{Kummer K3 surfaces.}\label{k.k3}
 Throughout this subsection $\mathrm{char} \, K \ne 2$.  Let $A$ be a  $2$-dimensional abelian variety over $K$. Then the group subscheme $A_2:=\Ker(A\rar{[2]}A)\subset A$ of 2-torsion points is reduced of order $16$.
The quotient $A/\sigma$  modulo the involution $A\rar{\sigma}A$, $\sigma(x) =-x$, is a projective surface, whose singular locus is precisely the image of $A_2$.   A  Kummer K3 surface $X$ is the blow up of $A/\sigma$ at
 $A_2 \mono  A/\sigma$ which is smooth.  Any translation invariant  differential $2$-form on $A$ descends to a non-vanishing regular  form $\omega$ on $X$.
Equivalently, $X$ can be viewed as the quotient of the variety $Z$ obtained  from $A$ by blowing up at $A_2$.
\begin{Th}\label{kum.k3.surf}
 Conjecture \ref{mainc} is true if $X$ is a Kummer K3 surface and $\mathrm{char} \, k \ne 2$. 
\end{Th}
\begin{proof}
%Fix a $2$-dimensional abelian  variety $A$ the surface $X$ is associated with and 
Fix a prime number $\ell \ne \mathrm{char}\, k$. By Theorem \ref{main2} since $X$ is maximally degenerate and, for some finite extension $K'\supset K$, the $\mathrm{Gal}(\oK/K')$-module $H^2(X_\oK, \bQ_{\ell})$ is isomorphic to $H^2(A_\oK, \bQ_{\ell})\oplus \bQ_{\ell}(-1)^{\oplus 16}$
the abelian variety $A$ is maximally degenerate.
Thus, after replacing $K$ by its finite extension we may assume that the analytic space, $A^{an}$, is the quotient of a split $2$-dimensional torus $T^{an}$ by a split lattice $\Lambda \subset T(K)$.  We also assume that all  the $2$-torsion points of $A$ are defined over $K$.
Under these assumptions we will prove that the formula (\ref{conjfor}) is true for $L=K$ and therefore for all its finite extensions. To do this we shall construct a formal poly-stable model $\oX$ for the analytic space $X^{an}$. By a general result of Berkovich
(\cite{ber3}, \S 5) the topological space $|X^{an}_\sK|$ is homotopy equivalent to the realization of the nerve of the special fiber\footnote{The notion of nerve of a scheme is recalled in Remark \ref{nerve}.} $\oX\times \overline k$, denoted by $|N(\oX\times \overline k)|$. On the other hand, the smooth locus of such a model is a weak N\'eron model of  $X^{an}$ and, thus, can be used to compute the motivic integral.

Let $\Xi$ be the group of characters of $T$ and $\Xi ^*$ the dual group. We have a canonical injective homomorphism $\rho:  \Lambda \to \Xi^*$ given by the valuation on $K$. Choose bases $\{v_1, v_2\}$,  $\{u_1, u_2\}$ for  $\Lambda$ and $\Xi$ such that the matrix of $\rho$ is diagonal
\[ \left( \begin{array}{cc}
m_1& 0  \\
0 & m_2 \end{array} \right)\] 
with positive $m_i$,  and let 
$T\simeq \bG_{m, K}\times \bG_{m, K}$, $\Lambda \simeq \bZ^2$ be the corresponding isomorphisms. Consider the standard  ``relatively complete'' model $\overline \bG_m$ of $ \bG_{m, K}$ over $R$ in the sense of Mumford (\cite{mu}, \S 5). We view $\overline \bG_m$ as a formal scheme over $R$ whose associated $K$-analytic space
is $\bG_{m, K}^{an}$ and which is equipped with an action of the multiplicative group $K^*$, extending the translation action on $\bG_{m, K}^{an}$, and an involution
that acts as $x\mapsto x^{-1}$ on $\bG_{m, K}^{an}$.
The  special fiber of $\overline \bG_m$ is a chain of projective lines $\bP^1_k$;  the action of $K^*$ induces a simple and transitive
action of $\bZ\simeq K^*/R^*$ on the set of its irreducible components. 
The quotient of   $\overline \bG_m \times  \overline \bG_m$ by $\Lambda \subset T(K)=K^*\times K^* $ is a proper strictly poly-stable formal model of $A^{an}$. Its smooth locus is the formal N\'eron model of $A^{an}$ (\cite{bs},  Def. 1.1 and Theorem 6.2). In particular, our assumptions on $A$ and $K$ imply that the $2$-torsion points of $A$
define $16$ sections of $(\overline \bG_m \times  \overline \bG_m)/ \Lambda$ over $R$ meeting the special fiber at distinct smooth points. Let $\oZ$ be the blow up of   $(\overline \bG_m \times  \overline \bG_m)/ \Lambda$ at these sections.  
By construction, the involution $\sigma$ of
$A$ extends to an involution $\overline \sigma$ of $\oZ$. The quotient  $\oX= \oZ/\overline \sigma $ is again a proper strictly poly-stable formal scheme whose generic fiber is the analytic K3 surface $X^{an}$. 
In particular, every $K'$-point  of $X^{an}$, where $K'\supset K$ is a finite unramified extension of $K$ reduces to a nonsingular point of the special fiber of $\oX$.  It follows that the smooth locus $\oX_{sm}\subset \oX$ of $\oX$ is a formal weak N\'eron model of $X^{an}$ (\cite{bs},  Def. 1.3).
As the N\'eron top degree differential form on $A$ induces a regular non-vanishing
differential form on $\oX_{sm}$ using (\cite{ls},  Theorem 4.4.1) we see that the motivic integral $\int_{X}$ equals the class $[Y_{sm}]$ of the smooth locus of the special fiber $Y$ of $\oX$.  We shall show  that $\int_{X}$ depends only on the order of the group $C=\bZ/m_1\bZ \times \bZ/m_2\bZ $ 
of connected components of the formal  N\'eron model of $A^{an}$. Indeed, since  all the 2-torsion points of $A$ are defined over $K$, the numbers $m_i$ are even. Thus, the involution $\sigma$ has precisely $4$  fixed points on $C$. It follows that  $Y_{sm}$ has $\frac{|C|}{2}+2$ connected components. All the components of $Y_{sm}$ are isomorphic to $\bG_{m, k}\times \bG_{m, k}$ except for those $4$ that correspond to fixed points of $\sigma$ on $C$. These $4$ components are isomorphic to the blow up of $\bG_{m, k}\times \bG_{m, k}$ at $4$ points of order $2$. Summarizing, we find
 $$\int_{X }= \left(\frac{|C|}{2}+2 \right)\bQ(0)+(20- |C| )\bQ(-1)+\left( \frac{|C|}{2}+2  \right)\bQ(-2).$$
Thus, to complete the proof of the formula (\ref{conjfor}) we need to show that $|C|=r_2(X,K)$.
Consider the commutative diagram induced by the morphism $f: \oZ \to \oX $ of the formal schemes 
    \begin{equation}\label{pks}
\def\normalbaselines{\baselineskip20pt
\lineskip3pt  \lineskiplimit3pt}
\def\mapright#1{\smash{
\mathop{\to}\limits^{#1}}}
\def\mapdown#1{\Big\downarrow\rlap
{$\vcenter{\hbox{$\scriptstyle#1$}}$}}
\begin{matrix}
H^2 (|N(\oX \times \overline k)|, \bZ )&   \iso &  \Gamma^2_\bZ(X) & \rar {}& H^2( X^{an}_\sK, \bZ_{\ell}) &\simeq &  H^2( X_\oK, \bZ_{\ell})\cr
  \mapdown{f^* }& & \mapdown{f^* } & & \mapdown{f^* } & & \mapdown{f^* } \cr
 H^2_{Zar}(|N(\oZ \times \overline k)|, \bZ )  & \iso & \Gamma^2_\bZ(Z) & \rar{ } &  H^2( Z^{an}_\sK, \bZ_{\ell})&\simeq &  H^2( Z_\oK, \bZ_{\ell}).
\end{matrix}
 \end{equation}
The topological space $|N(\oZ \times \overline k)|$ is homeomorphic to a real $2$-dimensional torus\footnote{Indeed, the scheme $\oZ \times \overline k$ is isomorphic to a direct product
of two (reducible) curves $D_i$, $i=1,2$, which, in turn, are isomorphic to $m_i$-gons of $\bP^1_{\overline k}$'s.  Thus, using that the formation $|N(\cdot)|$ commutes with products of
poly-stable schemes over algebraically closed field (\cite{ber3},  Prop. 3.14 (ii) and Cor. 3.17), we find that
$|N(\oZ \times \overline k)|\simeq |N(D_1)|\times |N(D_2)|\simeq S^1\times S^1$.};  the map $ |N(  \oZ \times \overline k)|\to |N(  \oX \times \overline k)|$ induced by $f$ identifies the target space with the quotient of the torus modulo the involution that takes a point to its inverse (with respect to the usual group structure on the real torus).
In particular,  $|N(  \oX \times \overline k)|$ is homeomorphic to a $2$-dimensional sphere. This proves the first part of the Theorem. Moreover, we have a commutative diagram
    \begin{equation}\label{compmonpar}
\def\normalbaselines{\baselineskip20pt
\lineskip3pt  \lineskiplimit3pt}
\def\mapright#1{\smash{
\mathop{\to}\limits^{#1}}}
\def\mapdown#1{\Big\downarrow\rlap
{$\vcenter{\hbox{$\scriptstyle#1$}}$}}
\begin{matrix}
\bZ &   \iso &  \Gamma^2_\bZ(X)&& \cr
  \mapdown{2 }& & \mapdown{f^* }&&  \cr
 \bZ  & \iso & \Gamma^2_\bZ(Z)&\simeq & \Gamma^2_\bZ(A),
\end{matrix}
 \end{equation}
 where the isomorphism $\Gamma^2_\bZ(A) \iso \Gamma^2_\bZ(Z)$ is induced by the morphism of formal schemes $\oZ \to (\overline \bG_m \times  \overline \bG_m)/ \Lambda $ that identifies the nerves of their special fibers.
On the other hand, since the morphism $Z \to X$ induced by $f$ has degree $2$, we have
$$(x,y)= \frac{1}{2}(f^*(x), f^*(y)), \quad x,y \in \Gamma^2_\bZ(X).$$
Comparing this  with (\ref{compmonpar}) we find that
$$r_2(X, K)=\frac{r_2(Z, K)}{2} = \frac{r_2(A, K)}{2}.$$
It remains to show that $\frac{r_2(A, K)}{2}= |C|$. Consider the exact sequence of $G$-modules
$$ 0\to \Lambda ^* \otimes \bZ_{\ell} \to H^1(A_\oK, \bZ_{\ell}) \to \Xi  \otimes \bZ_{\ell}(-1) \to 0.$$
The canonical morphism $\Gamma^1_\bZ(A) \mono H^1(A_\oK, \bZ_{\ell})$ identifies $\Gamma^1_\bZ(A)$ with $\Lambda ^*\subset H^1(A_\oK, \bZ_{\ell})$ (\cite{ber1}, \S 6.5).
Let $\tilde u_1 \wedge \tilde u_2$ be an element of $H^2(A_\oK, \bZ_{\ell})(2)$ that projects to $  u_1\wedge u_2 \in \bigwedge^2\Xi \otimes \bZ_{\ell}$. Then, we have
$$ N^2(\tilde u_1 \wedge \tilde u_2)= N(N(\tilde u_1)\wedge \tilde u_2 + \tilde u_1 \wedge N(\tilde u_2))= 2 N(\tilde u_1) \wedge N(\tilde u_2)=$$
$$ 2m_1m_2 (v^*_1 \wedge  v^*_2) \in \bigwedge^2\Lambda ^* \otimes \bZ_{\ell}.$$
It follows that the monodromy pairing on  $\Gamma^2_\bZ(A)\simeq \bigwedge^2\Lambda ^*$ is given by the formula
$$(v^*_1\wedge v^*_2, v^*_1\wedge v^*_2)=  -\frac{<v^*_1\wedge v^*_2, \tilde u_1 \wedge \tilde u_2>}{2m_1m_2}=  \frac{1}{2m_1m_2}$$
and therefore $r_2(A, K)= 2m_1m_2= 2|C|$. 
 The proof of the formula (\ref{conjfor}) is completed. 
 
 Let us prove the second statement of the theorem. We will derive it from an analogous result for abelian varieties  (\cite{ber1}, \S 6.5) which asserts that
 the lattice  
 $$ \Gamma^2_{\bZ_{\ell}}(A) \mono H^2( A_{\overline K}, \bZ_{\ell})$$
 is saturated.
 It follows that the lattice
 $$ \Gamma^2_{\bZ_{\ell}}(A)  \iso \Gamma^2_{\bZ_l}(Z)    \mono H^2( Z_\oK, \bZ_{\ell})\simeq  H^2( A_\oK, \bZ_{\ell}) \oplus \bZ_{\ell}(-1)^{\oplus 16}$$
 is also saturated.
 We claim that in the commutative diagram
     \begin{equation}\label{repeat}
\def\normalbaselines{\baselineskip20pt
\lineskip3pt  \lineskiplimit3pt}
\def\mapright#1{\smash{
\mathop{\to}\limits^{#1}}}
\def\mapdown#1{\Big\downarrow\rlap
{$\vcenter{\hbox{$\scriptstyle#1$}}$}}
\begin{matrix}
\bZ_{\ell}&   \iso &  \Gamma^2_{\bZ_{\ell}}(X) & \rar {}&  H^2( X_\oK, \bZ_{\ell})\cr
  \mapdown{2 }& & \mapdown{f^* } &  & \mapdown{f^* } \cr
 \bZ_{\ell}  & \iso & \Gamma^2_{\bZ_{\ell}}(Z)  & \rar{ } &   H^2( Z_\oK, \bZ_{\ell})   .
\end{matrix}
 \end{equation}
the vertical morphisms are isomorphisms up to $2$-torsion. Indeed, the compositions $f_* f^*$ and $f^* f_*$ with the trace morphism $H^2( Z_\oK, \bZ_{\ell}) \rar{f_*}  H^2( X_\oK, \bZ_{\ell})$
are equal to $2\, Id$. This proves the second part of the theorem for $\ell\ne 2$.
For $\ell=2$, we apply a result of Nikulin  (see, {\it e.g.} \cite{mo2}, Lemma on p. 56) that states the lattice
   \begin{equation}\label{cohofkum}
H^2( A_\oK, \bZ_{\ell}) \mono H^2( Z_\oK, \bZ_{\ell})  \stackrel{f_*}{\mono} H^2( X_\oK, \bZ_{\ell}) 
  \end{equation}
  is saturated. As the lattice $\Gamma^2_{\bZ_{\ell}}(X)$ equals the image of the saturated lattice $\Gamma^2_{\bZ_{\ell}}(A)  \mono H^2( A_{\overline K}, \bZ_{\ell})$ under composition (\ref{cohofkum}) it is saturated as well.
\end{proof}

\newpage
\appendix
\section{Erratum}
\begin{abstract}
It was pointed out to us by Olivier Wittenberg that Kulikov's Theorem \cite{ku}, which we used in the proof of two results,  is stated in \cite{sv} incorrectly.  We give a correct statement of the Kulikov Theorem  and show that both results remain valid. 

 \end{abstract}
%\section{Introduction}
Let $\oC$ be smooth curve over $\bC$, and let  $a\in \oC$ be a closed point, $C=\oC-a$.  Denote by $R$ the completion of the local ring $ \cO_{\oC, a}$ and by $K$ its field of fractions.
In \S 3.2 of \cite{sv} we made the following assertion referring to \cite{ku}:
\begin{assertion}\label{a1} Let $\overline \cX \to \oC$ be a  projective strictly semi-stable morphism whose restriction $\cX$ to $C$ is a smooth family of K3 surfaces. Then, there is 
 a {\it projective} strictly semi-stable family $\overline \cX' \to \oC$, whose restriction to $C$ is isomorphic to $\cX$, and the log canonical bundle  $\Omega^2_{\overline \cX'/\oC}(log)$ is trivial over a neighborhood of the special fiber 
 of  $\overline \cX' $.
 \end{assertion}
 As it was pointed out to us by Olivier Wittenberg, the actual result proven \cite{ku}  is weaker then the above assertion. A correct statement of Kulikov's theorem reads as follows:
 \begin{theo}[Kulikov]\label{kulikov}  Let $\overline \cX \to \oC$ be a  projective strictly semi-stable morphism whose restriction $\cX$ to $C$ is a smooth family of K3 surfaces.
 Then there is a proper strictly semi-stable  complex analytic space  over $\oC$, $\overline \cX' \to \oC$, together with a bimeromorphic map
 $\overline \cX' \dashrightarrow \overline \cX$, which commutes with the projections to $\oC$ and induces an isomorphism  $\overline \cX' \times_\oC C \iso \cX \times_\oC C$, such that 
 the log canonical bundle  $\Omega^2_{\overline \cX'/\oC}(log)$ is trivial over a neighborhood of the special fiber 
 of  $\overline \cX' $. 
\end{theo}
One refers to $\overline \cX' $  as a Kulikov model for $\cX \to C$.  It is shown in \cite{ku} that the special fiber, $Y'$, of any Kulikov model has a very special form (this part of Kulikov's Theorem is stated correctly in  \cite{sv},  \S 3.2).
In particular,  all the irreducible components of $Y'$ are smooth projective surfaces.

We do not know whether Assertion \ref{a1} is true.  In fact, we do not even know whether, for given $\cX \to C$, there exists a Kulikov model which is a scheme. 

However, we are going to show that Theorem 1 and Proposition 4.8 from \cite{sv}, which are the only results of {\it loc. cit.} whose proofs were based on Assertion \ref{a1}, are still valid and can be proved using
 Kulikov's Theorem. In fact, an inspection of the proof of Theorem 1  given in \cite{sv}  shows  that  only part of the argument where algebraicity of Kulikov's model $\overline \cX'$ is used is the proof of formula (3.3) stated below.
  \begin{lm}\label{lemma1er} If $\overline \cX'$ is a Kulikov model and $Y'$ is its special fiber, we have
 \begin{equation}\label{minkul}
 \int_{X}=[Y'_{sm}].
  \end{equation}
  Here $X$ is the K3 surface over $K$  obtained from $\cX$ by the base change.
 \end{lm}
Note that since the irreducible components of $Y'$ are projective, the smooth locus of $Y'_{sm}$ has the structure of an algebraic variety and, hence, $[Y'_{sm}]$ makes sense as an element of the Grothendieck  group of varieties.
 \begin{proof}[Proof of Lemma \ref{lemma1er}]
 Let $\overline \cX' \to \oC$ be any proper  complex analytic model for $\cX \to C$, which is bimeromorphically isomorphic to an algebraic model $\overline \cX \to \oC$.  Then, according to a theorem of Artin  (\cite{ar}, Theorem 7.3),
 $\overline \cX'$ has a unique structure of an algebraic space such that  the map $\overline \cX' \to \oC$ is algebraic. Since any smooth algebraic surface is quasi-projective, the smooth locus $Y_{sm}'$ of
 the special fiber of  $\overline \cX' \to \oC$ acquires the structure of an algebraic variety. Now, assume that $\overline \cX'$ is regular and let $V_i^\circ$ be irreducible components of  $Y'_{sm}$. Define integers $m_i$  by  the formula
 (1.1) from \cite{sv}. The lemma will follow from a more general claim:
 \begin{equation}\label{n.m.i.e}
\int_{X}: =  \sum _i   [ V_i^\circ  ] ( m_i  - \min_i  m_i ). 
\end{equation}
To prove (\ref{n.m.i.e}) observe, first, that the formula is true if $\overline \cX' $ is a scheme. Indeed, in this case $\overline \cX'  \otimes R -Y'_{sing}$ is a weak N\'eron model  (\cite{sv}, Lemma 2.10), and (\ref{n.m.i.e}) boils down to the definition of motivic integral 
(\cite{sv}, \S 1.1). 
 Hence, it suffices to check that the right-hand side of (\ref{n.m.i.e})  is independent of the choice of a regular model $\overline \cX' $.
Using the Weak Factorization Theorem for algebraic spaces (\cite{akmw}, Theorem 0.3.1) it is enough to show this for two models, one of which is obtained from the other one by blowing up at a smooth subvariety of the special fiber. 
In this case the assertion follows by a direct inspection. 
  \end{proof}
 Next, we correct the proof of Proposition 4.8  given in  \cite{sv}. It suffices to prove the following result.
 \begin{lm}\label{lemma2er}
  With the assumption of Kulikov's Theorem,  $X:= \cX \times _C \text{ spec}\, K $, there is a homotopy equivalence between the topological space  $|X^{an}_\sK |$  and the Clemens polytope $Cl(Y')$ of the special fiber of $\overline \cX'$, which
  identifies the canonical map 
$$\Gamma^*_\bZ(X) \to  H^*(\lim X,\bZ) $$
with the map $$H^*(Cl(Y')) \to H^*(\lim X,\bZ)$$ 
coming from the weight spectral sequence. 
\end{lm}
 \begin{proof} Berkovich's result (\cite{ber3}, \S 5) applied to the strictly semi-stable model $\overline \cX $ (which is a scheme) implies 
 the assertion of the lemma  with $Y'$ replaced by $Y$.
 Next, let  $\overline \cX\to \oC$ and $ \overline \cX' \to \oC$ be any proper regular models for $\cX \to C$ in the category of algebraic spaces.  Assume that  the {\it reduced} special fibers $Y_{red}$ and $Y'_{red}$ are strict normal crossing divisors. We claim that there
is a homotopy equivalence between
the Clemens polytopes   $Cl(Y)$ and $Cl(Y')$, which identifies the maps  $H^*(Cl(Y))\to H^*(\lim X,\bZ), H^*(Cl(Y')) \to H^*(\lim X,\bZ)$.
Indeed, using the Weak Factorization Theorem  (\cite{akmw}, Theorem 0.3.1) we may assume that $\overline \cX$ is obtained from $\overline \cX'$
by  blowing up at an admissible subvariety of the special fiber. In this case the assertion can be checked directly (see \cite{s}, Lemma in \S 2).  
\end{proof}
\begin{rem}\label{rem.e}  \begin{enumerate}[(a)]
 \item The first part of Lemma \ref{lemma2er}  is a corollary of an unpublished result of Michael Temkin:  if $X$ is a smooth proper scheme over the fraction field $K$ of  complete discrete valuation ring $R$ with perfect residue field 
$k$ and $\oX$ is a proper  strictly semi-stable algebraic space  over $R$ with the generic fiber isomorphic to $X$,  then the space  $|X^{an}_\sK |$ is homotopy equivalent to the Clemens polytope $Cl(Y_{\ok})$ of the geometric special fiber.
\item  We expect that the motivic integral of a  Calabi-Yau variety  can be computed from its weak N\'eron model in the category of algebraic spaces: 
 if $X$ is a smooth proper Calabi-Yau variety over $K$, $\cV$ is an algebraic space over $R$ which is a weak N\'eron model for $X$, and $V_i^\circ $ are irreducible components of the special fiber of $\cV$, then formula
  (\ref{n.m.i.e}) is true. Note that 
the natural homomorphism from the Grothendieck group of varieties to the
 Grothendieck group of algebraic spaces of finite type over a field is isomorphism ({\it cf.} \cite{kn}, Proposition 2.6.7).  
   \end{enumerate}
\end{rem}  
  
   \section*{Acknowledgements}
 We are grateful to Olivier Wittenberg for pointing out our misinterpretation of Kulikov's Theorem and for his suggestion to use algebraic spaces to overcome the araised problem. We are also grateful to Michael Temkin
 for explaining to us his result stated in Remark \ref{rem.e}.
 
 \end{document}